%% file: Integral_reps.tex
\documentclass{article}
\usepackage{array}
\usepackage{amsmath}
\usepackage{amssymb}
\usepackage{graphicx}
\usepackage{url}
\usepackage{maplestd2e}

\DefineParaStyle{Maple Heading 1}
\DefineParaStyle{Maple Text Output}
\DefineParaStyle{Maple Dash Item}
\DefineParaStyle{Maple Bullet Item}
\DefineParaStyle{Maple Normal}
\DefineParaStyle{Maple Heading 4}
\DefineParaStyle{Maple Heading 3}
\DefineParaStyle{Maple Heading 2}
\DefineParaStyle{Maple Warning}
\DefineParaStyle{Maple Title}
\DefineParaStyle{Maple Error}
\DefineCharStyle{Maple Hyperlink}
\DefineCharStyle{Maple 2D Math}
\DefineCharStyle{Maple Maple Input}
\DefineCharStyle{Maple 2D Output}
\DefineCharStyle{Maple 2D Input}
\usepackage[top=.9in,bottom=.9in]{geometry}
 \numberwithin{equation}{section}
 
\DefineParaStyle{Maple Output}
\DefineCharStyle{2D Comment}
\DefineCharStyle{2D Math}
\DefineCharStyle{2D Output}
\begin{document}
\pagestyle{plain}
\begin{maplegroup}
\widowpenalty=600
\clubpenalty=600
\title{Integral and Series Representations of Riemann's Zeta function,  Dirichelet's Eta Function and a Medley of Related Results}

\author{Michael S. Milgram \footnote{mike@geometrics-unlimited.com} \\
\it Consulting Physicist, Geometrics Unlimited, Ltd. \\
\it Box 1484, Deep River, Ont. Canada. K0J 1P0 }
\date{August 13, 2012}
\maketitle

\setcounter{page}{1}
\pagenumbering{arabic}

\begin{flushleft}
\begin{abstract}

Abstract: Contour integral representations for Riemann's Zeta function and Dirichelet's Eta (alternating Zeta) function are presented and investigated. These representations flow naturally from methods developed in the 1800's, but somehow they do not appear in the standard reference summaries, textbooks or literature. Using these representations as a basis, alternate derivations of known series and integral representations for the Zeta and Eta function are obtained on a unified basis that differs from the textbook approach, and results are developed that appear to be new. 

\end{abstract}

\section{Introduction}

Riemann's Zeta function $\zeta(s)$, and its sibling, Dirichelet's (alternating zeta) function $\eta(s)$ play an important role in physics, complex analysis and number theory, and have been studied extensively for several centuries. In the same vein, the general importance of a contour integral representation for any function has been known for almost two centuries, so it is surprising that contour integral representations for both $\zeta(s)$ and $\eta(s)$ exist that cannot be found in any of the modern handbooks (NIST, \cite{OLBC}, Section 25.5(iii); Abramowitz and Stegun, \cite{Abram}, Chapter 23), textbooks (Apostol, \cite{Apostol}, Chapter 12; Olver, \cite{Olver}, Chapter 8.2; Titchmarsh, \cite{Titchmarsh}, Chapter 4; Whittaker and Watson,  \cite{WhitWat}, Section (13.13)), summaries (Edwards, \cite{Edwards}, Ivic, \cite{Ivic}, Chapter 4; Patterson, \cite{Patterson}; Srivastava and Choi, \cite{Sriv} and \cite{SrivChoi2}), compendia (Bateman, \cite{Bateman}, Section 1.12 and Chapter 17), tables (Gradshteyn and Ryzhik, \cite{G&R}, 9.512; Prudnikov et. al, \cite{Prudnikov}, Appendix II.7) and websites (\cite{Encyclop}, \cite{OLBC},   \cite{WikiZ}, \cite{WikiN}, \cite{Hankel}, \cite{Wolfram}) that summarize what is known about these functions. One can only conclude that such representations, generally discussed in the literature of the late 1800's, have been long-buried and their significance overlooked by modern scholars.\newline 

It is the purpose of this work to disinter these representations, revisit and explore some of the consequences. In particular, it is possible to obtain series and integral representations of $\zeta(s)$ and $\eta(s)$ that unify and generalize well-known results in various limits and combinations and reproduce results that are only now being discovered and investigated by alternate means. Additionally, it is possible to explore the properties of $\zeta(s)$ in the complex plane and on the critical line, and obtain results that appear to be new. Many of the results obtained are disparate and difficult to categorize in a unified manner, but share the common theme that they are all somehow obtained from a study of the revived integral representations. That is the unifying theme of this work. To maintain a semblance of brevity, most of the derivations are only sketched, with citations sufficient enough to allow the reader to reproduce any result for her/himself.

\input{contour_rep}

\section{Simple reductions}
\input{simple_0_leq_c_leq_1.tex}

\subsection{Simple Recursion via Partial Differentiation}
\input{more_simple.tex}

\section{Special cases}
\subsection{The case c=0}
\input{section_c=0.tex}

\subsection{The case $c\geq 1$}\label{sec:c geq 1}
\input{c_geq_1.tex}

\section{Series Representations}
\subsection{By Parts}
\input{by_parts.tex}

\subsection{By splitting}\label{sec:splitting}
\input{split2.tex}

\section{Integration by parts}
\input{by_parts_2.tex}

\section{Conversion to Integral Representations}
\input{int_reps.tex}

\section{The case $s=n$} \label{sec:s=n}
\input{s_eq_n.tex}

\section{The Critical Strip}\label{sec:crits}
\input{critstrip.tex}

\section{Comments}
\input{comments.tex}

\appendix
\section{Appendix}
\input{appendix.tex}

\end{flushleft} 
\end{maplegroup}
\end{document}

%% file: contour_rep.tex
\section{Contour integral representations for $\zeta(s)$ and $\eta(s)$ }

Throughout, I shall use $s=\sigma + i\rho$, where $s \in \mathfrak{C}$ (complex), $\sigma , \rho \in \mathfrak{R}$ (reals) and $n$ and $N$ are positive integers.\newline

The Riemann Zeta function $\zeta(s)$ and the Dirichelet alternating zeta function  $\eta(s)$ are well-known and defined by (convergent) series representations:

\begin{equation}
\zeta(s) \equiv \sum_{k=1}^{\infty} {k^{-s}} \hspace*{7 cm}
    {\sigma > 1}
\label{ZetaSum}
\end{equation}

\begin{equation}
\eta(s) \equiv -\sum_{k=1}^{\infty} {(-1)^k\,k^{-s}} \hspace*{6cm} {\sigma > 0}
\label{EtaSum}
\end{equation}

with
\begin{equation}
\eta(s) = (1-2^{1-s})\,\zeta(s)
\label{Eta=Zeta} \,.
\end{equation} 

It is easily found by an elementary application of the residue theorem, that the following reproduces (\ref{ZetaSum})

\begin{equation}
\zeta(s) = \pi \sum_{k=1}^{\infty} residue \vert_{t=k} \oint_{t \circeq k} \frac{t^{-s} e\,^{i\pi t}}{\mathrm{sin}(\pi t)}\, dt \hspace*{4cm}
    {\sigma > 1} 
\label{SumRes1}
\end{equation}
and the following reproduces (\ref{EtaSum})
\begin{equation}
\eta(s) = \pi \sum_{k=1}^{\infty} residue \vert_{t=k} \oint_{t \circeq k} {t^{-s} \,\mathrm{cot}(\pi t)\,e\,^{i\pi t}} dt \hspace*{3 cm} {\sigma > 0} 
\label{SumAltRes1}
\end{equation}

where the contour of each of the integrals encloses the pole at $t=k$.\newline

Convert the result (\ref{SumRes1}) into a contour integral enclosing the positive integers on the t-axis in a clockwise direction, and, provided that $\sigma > 1$ so that contributions from infinity vanish, the contour may then be opened such that it stretches vertically in the complex t-plane, giving

\begin{eqnarray}
{\zeta(s) = \frac{i}{2} \int_{c-i\infty}^{c+i\infty}\,t^{-s}\,\mathrm{cot}(\pi t)\, dt} \hspace*{4 cm} 0<c<1, c\in \mathfrak{R} \label{ContInt0} \\ 
{= -\frac{1}{2} \int_{-\infty}^{\infty}\,(c+it)^{-s}\,\mathrm{cot}(\pi(c+i t))\, dt} \hspace*{4 cm} \sigma > 1 ,
\label{ContInt}
\end{eqnarray}

a result that cannot be found in any of the modern reference works cited (in particular \cite{SG&A}). It should be noted that this procedure is hardly novel - it is often employed in an educational context (e.g. in the special case $s=2$) to demonstrate the utility of complex integration to evaluate special sums (\cite{Sriv} Section 4.1), (\cite{Flaj&Salvy}, Lemma 2.1), and forms the basis for an entire branch of physics (\cite{Newton}). \newline

A similar application to (\ref{SumAltRes1}), with some trivial trigonometric simplification, gives
\begin{eqnarray}
{ \eta (s)={\displaystyle \frac {i}{2}} \,{\displaystyle \int _{c
 - i\infty }^{c + i\infty }} {\displaystyle \frac {t^{-s} }{
\mathrm{sin}(\pi \,t)}} \,dt }
 \label{ContIntEta0} \\ 
{\, \,\, ={\displaystyle \frac {1}{2}} \,{\displaystyle \int _{ - 
\infty }^{\infty }} {\displaystyle \frac {(c + i\,t)^{-s}}{
\mathrm{sin}(\pi \,(c + i\,t))}} \,dt \, \, .
\label{ContIntEta1}
}
\end{eqnarray}

Eq. (\ref{ContIntEta1}) reduces, with $c=1/2$, to a Jensen result (\cite{Jensen}), presented in 1895, and to the best of my knowledge, reproduced, again, but only in the special case $c=1/2$, only once in the modern reference literature (\cite{WikiN}). In particular, although Lindelof's 1905 work \cite{Lindelof} arrives at a result equivalent to that obtained here, attributing the technique to Cauchy (1827), it is performed in a general context, so that the specific results (\ref{ContInt}) and (\ref{ContIntEta1}) never appear, particularly in that section of the work devoted to $\zeta(s)$. In modern notation, for a meromorphic function $f(z)$, Lindelof writes (\cite{Lindelof}, Eq. {\it{III}} (4)) \textbf{in general}
 
\begin{equation}
\sum_{\nu =-\infty}^{\infty} f\,(\nu) = -\frac {1}{2 \, \pi \, i} \oint \pi \, \mathrm {cot}( \pi \, z )\, f(z) \, dz,
\label{LindEta} 
\end{equation}

where the contour of integration encloses all the singularities of the integrand, but never specifically identifies $f(z)=z^{-s}$. Much later, after the integral has been split into two by the invocation of identities for $cot(\pi\,z)$, he does make this identification (\cite{Lindelof}, Eq. {\it{4.III}} (6)), arriving (with \ref{Eta=Zeta}) at (\ref{EtaC=half} - see below) a now-well-known result, thereby bypassing the contour integral representation (\ref{ContInt}). In a similar vein, Olver (\cite{Olver}, page 290) reproduces (\ref{LindEta}), employing it to evaluate remainder terms, the Abel-Plana formula and a Jensen result, but again, never writes (\ref{ContInt}) explicitly. Perhaps the omission of an explicit statement of (\ref{ContInt}) and (\ref{ContIntEta1}) from Chapter IV of Lindelof (\cite{Lindelof}) is why these forms seem to have vanished from the historical record, although citations to Jensen \cite{Jensen} and Lindelof \cite{Lindelof} abound in modern summaries\footnote{Recently, Srivastava and Choi \cite{SrivChoi2}, pages 169-172 have reproduced Lindelof's analysis, also failing to arrive explicitly at (\ref{ContInt}); although they do write a general form of (\ref{ContIntEta0}) (\cite{SrivChoi2}, Section 2.3, Eq. (38)), it is studied only for the case $c=1/2.$}. \newline
 
At this point, it is worthwhile to quote several related integral representations that define some of the relevant and important functions of complex analysis, starting with Riemann's original 1859 contour integral representation of the zeta function (\cite{Riemann})
\begin{equation}
\zeta(s) =  \frac {\Gamma(1-s)}{2 \pi i} \oint \frac{t^{s-1}}{e^{-t}-1} {dt}
\label{ContRie}
\end{equation}

where the contour of integration  encloses the negative t-axis, looping from $t=-\infty - i0$  to $t=-\infty + i0$ enclosing the point t=0.  \newline

The equivalence of (\ref{ContRie}) and (\ref{ZetaSum}) is well-described in texts (e.g. \cite{SrivChoi2}) and is obtained by reducing three components of the contour in (\ref{ContRie}). A different analysis is possible however (e.g. \cite{Edwards}, Section 1.6) - open and translate the contour in (\ref{ContRie}) such that it lies vertically to the right of the origin in the complex t-plane\footnote{and thereby vanishes} with $\sigma < 0$ and evaluate the residues of each of the poles of the integrand lying at $t= \pm \, 2n\pi i$ to find

\begin{equation}
\zeta(s) = \frac{i}{2\pi}\Gamma(1-s)(2\pi i)^{s}(\exp (i\pi s/2)- \exp (-i\pi s/2))\sum_{k=1}^{\infty} k^{s-1}  .
\label{pre-functional}
\end{equation}

This can be further reduced, giving the well-known functional equation for $\zeta(s)$

\begin{equation}
\zeta(1-s) = \frac {2 \, \mathrm{cos}(\pi s /2) \, \Gamma(s)} {(2 \pi)^{s}} \, \zeta(s)
\label{functional}
\end{equation}

valid for all $s$ by analytic continuation from $\sigma < 0$. This demonstrates that (\ref{ContRie}) and (\ref{ContInt0}) are fundamentally different representations of $\zeta (s)$, since further deformations of (\ref{ContRie}) do not lead to (\ref{ContInt}).\newline

Ephemerally, it has been claimed  \cite{WikiN} that the following complex integral representation (missing a factor $2\pi$ in \cite{WikiN}) 

\begin{equation}
\zeta(s)=\frac{\pi}{2\,(s-1)} \int_{-\infty}^{\infty} \frac{(1/\,2+i\,t)^{1-s}}{\mathrm{cosh}^2(\pi\,t)}dt \hspace*{2cm} \sigma>0
\label{cosh^2_rep}
\end{equation}
can be obtained by the invocation of the Cauchy-Schl\"{o}milch transformation \cite{Amdeb}. In fact, the representation \eqref{cosh^2_rep} follows by the simple expedient of integrating \eqref{ContInt} by parts in the special case $c=1/2$ (see \eqref{Zetas+1} below), as also noted in \cite{WikiN}. Another result 

\begin{equation}
\zeta(s)=\frac{2^{s-1}}{(1-2^{1-s})\,\Gamma(s+1)}\int_{0}^{\infty}\frac{t^s}{\mathrm{cosh}^2(t)} dt \hspace*{2cm} \sigma>-1
\end{equation}

obtained in \cite{Amdeb} on the basis of the Cauchy-Schl\"{o}milch transformation does not follow from \eqref{cosh^2_rep} as erroneously claimed in \cite{WikiN}.\newline

Based on (\ref{ContRie}), a number of other representations of $\zeta(s)$ are well-known. From the tables cited, the following are worth noting here:\newline
(\cite{Bateman}, Eq.1.12(6)) 
\begin{equation}
\zeta(s)= \frac {1} {2 \Gamma(s) (1 - 2 ^{-s})} \int _{0}^{\infty} \frac {t^{s-1}} {\mathrm{sinh}(t)}  dt \,   \hspace*{3cm} \sigma >1\,,
\label{Erd1.12.6}
\end{equation}
\newline
(\cite{Bateman}, Eq.1.12(4)) 
\begin{equation}
\zeta(s)= \frac {1} {\Gamma(s) } \int _{0}^{\infty} \frac {t^{s-1}} {\mathrm{exp}(t)-1}  dt \,   \hspace*{3cm} \sigma >1\,,
\label{Erd1.12.4}
\end{equation}
\newline
and (\cite{G&R}, Eq. 3.527(1))
\begin{equation}
\mapleinline{inert}{2d}{Zeta(s) = 1/4*(2*a)^(s+1)/GAMMA(s+1)*Int(1/sinh(a*t)^2*t^s,t = 0 ..
infinity);}{%
\[
\zeta (s)={\displaystyle \frac {1}{4}} \, \! 
{\displaystyle \frac {(2\,a)^{(s + 1)}}{\Gamma (s + 1)}} \,
{\displaystyle \int _{0}^{\infty }} {\displaystyle \frac {t^{s}}{
\mathrm{sinh}^{2} (a\,t) }} \,dt \! \,\, .  \hspace*{3cm} \sigma >1
\]
\label{GR527.1}
}
\end{equation}

By way of comparison, Laplace's representation (\cite{G&R}, Eq. 8.315(2))

\begin{equation}
\frac{1}{\Gamma(s)} = \frac{1}{2 \pi i} \int_{-\infty}^{\infty} \frac{e^{(c+it)}} {(c+it)^s} dt
\label{GamInv}
\end{equation}
defines the inverse Gamma function.

%% file: simple_0_leq_c_leq_1.tex
\subsection{The case $0\leq c \leq 1$}
The simplest reductions of (\ref{ContInt}) and (\ref{ContIntEta1}) are obtained by uniting the two halves of their range corresponding to $t<0$ and $t>0$. With the use of

\begin{equation}
\mapleinline{inert}{2d}{(-(c+tau*I)^s+(c-I*tau)^s)*I =
2*(c^2+t^2)^(1/2*s)*sin(s*arctan(t/c));}{%
\[
( - (c + i \, t )^{s} + (c - i \, t )^{s})\,i=2\,(c^{2} + t 
^{2})^{( {s/2})}\,\mathrm{sin}(s\,\mathrm{arctan}(
{\displaystyle {t/c }} ))
\]
}
\label{ArcTanSin}
\end{equation}

and

\begin{equation}
{%
(c + i\,t)^{s} + (c - i\, t)^{s}=2\,(c^{2} + t^{2})^{(s/2
)}\,\mathrm{cos}(s\,\mathrm{arctan}({\displaystyle {t/c}} 
))
}
\label{ArcTanCos}
\end{equation}

one finds (with $0<c<1$ and $0 \leq c \leq 1 $ respectively)
\begin{equation}
\mapleinline{inert}{2d}{Zeta(s) =
-Int((c^2+tau^2)^(-1/2*s)*(cos(s*arctan(tau/c))*cos(Pi*c)*sin(Pi*c)-si
nh(Pi*tau)*cosh(Pi*tau)*sin(s*arctan(tau/c)))/(cosh(Pi*tau)^2-cos(Pi*c
)^2),tau = 0 .. infinity);}{%
\maplemultiline{
\zeta (s)= - {\displaystyle \int _{0}^{\infty }} {\displaystyle
\frac {\,\mathrm{cos}(s\,
\mathrm{arctan}({\displaystyle {t/c }} ))\,\mathrm{cos}
(\pi \,c)\,\mathrm{sin}(\pi \,c) - \mathrm{sinh}(\pi \,t )\,
\mathrm{cosh}(\pi \,t )\,\mathrm{sin}(s\,\mathrm{arctan}(
{\displaystyle {t/c }} ))}{(c^{2} + t^{2})^{s/2}\,(\mathrm{cosh}^{2} (\pi \,t )
- \mathrm{cos}^{2}(\pi \,c) )}}dt  \\ \hspace*{12cm}{\sigma>1}
  }
}
\label{Zeta1WithC}
\end{equation}
and
\begin{equation}
\mapleinline{inert}{2d}{eta(s) =
Int((-sin(s*arctan(t/c))*cos(Pi*c)*sinh(Pi*t)+cos(s*arctan(t/c))*sin(P
i*c)*cosh(Pi*t))/((c^2+t^2)^(1/2*s))/(cosh(Pi*t)^2-cos(Pi*c)^2),t = 0
.. infinity);}{%
\[
\eta (s)={\displaystyle \int _{0}^{\infty }} {\displaystyle 
\frac { - \mathrm{sin}(s\,\mathrm{arctan}({\displaystyle t/c } ))\,\mathrm{cos}(\pi \,c)\,\mathrm{sinh}(\pi \,t) + 
\mathrm{cos}(s\,\mathrm{arctan}({\displaystyle t/c } ))\,
\mathrm{sin}(\pi \,c)\,\mathrm{cosh}(\pi \,t)}{(c^{2} + t^{2})^{
s/2 }\,(\mathrm{cosh}^{2} (\pi \,t) - \mathrm{cos}^{2} (\pi \,
c) )}} \,dt\; ,
\]
}
\label{Eta1WithC}
\end{equation}

the latter valid for all $s$. In the case that $c=\frac{1}{2}$, (\ref{Zeta1WithC}) reduces to
 

\begin{equation}
\mapleinline{inert}{2d}{Zeta(s) =
2^s*Int(sinh(Pi*tau)/((1+4*tau^2)^(1/2*s))/cosh(Pi*tau)*sin(s*arctan(2
*tau)),tau = 0 .. infinity);}{%
\[
\zeta (s)=2^{s}\,{\displaystyle \int _{0}^{\infty }} 
{\displaystyle \frac {\mathrm{sinh}(\pi \,t )\,\mathrm{sin}(s
\,\mathrm{arctan}(2\,t ))}{(1 + 4\,t ^{2})^{s/2}
\,\mathrm{cosh}(\pi \,t )}} \,dt , \hspace*{2 cm} \sigma > 1 
\]
}
\label{JensenZ}
\end{equation}

equivalent to the well-known result\footnote{attributed to Jensen by Lindelof, and a special case of Hermite's result for the Hurwitz zeta function (\cite{WhitWat}, Section 13.2).} (which also appears as an exercise in many textbooks - e.g. \cite{SrivChoi2}, Eq.2.3(41))

\begin{equation}
\mapleinline{inert}{2d}{Zeta(s) =
2^(s-1)/(s-1)-2^s*int(sin(s*arctan(2*t))*exp(-Pi*t)/((1+4*t^2)^(1/2*s)
)/cosh(Pi*t),t = 0 .. infinity);}{%
\[
\zeta (s)={\displaystyle \frac {2^{(s - 1)}}{s - 1}}  - 2^{s}\,
{\displaystyle \int _{0}^{\infty }} {\displaystyle \frac {
\mathrm{sin}(s\,\mathrm{arctan}(2\,t))\,e^{ - \pi \,t}}{(1 + 4
\,t^{2})^{s/2}\,\mathrm{cosh}(\pi \,t)}} \,dt \;.
\]
}
\label{Sr2.3.41}
\end{equation}

 Similarly, (\ref{Eta1WithC}) reduces to the known (\cite{Lindelof}, {\it{ IV.III(6)}}) result 

\begin{equation}
\mapleinline{inert}{2d}{eta(s) = 2^s*Int(cos(s*arctan(2*t))/((1+4*t^2)^(1/2*s))/cosh(Pi*t),t
= 0 .. infinity);}{%
\[
\eta (s)=2^{s-1}\,{\displaystyle \int _{0}^{\infty }} 
{\displaystyle \frac {\mathrm{cos}(s\,\mathrm{arctan}(t))}{(1
 + t^{2})^{s/2}\,\mathrm{cosh}(\pi \,t\,/2)}} \,dt \,\, .
\]
}
\label{EtaC=half}
\end{equation}

Comparison of (\ref{JensenZ}) and \eqref{Sr2.3.41} immediately yields the known result (\cite{G&R}, Eq.3.975(2))

\begin{equation}
\mapleinline{inert}{2d}{Zeta(s) =
2^(s-1)/(1-s)+2^s*int(sin(s*arctan(2*t))*exp(Pi*t)/((1+4*t^2)^(1/2*s))
/cosh(Pi*t),t = 0 .. infinity);}{%
\[
\zeta (s)={\displaystyle \frac {2^{(s - 1)}}{1 - s}}  + 2^{s}\,
{\displaystyle \int _{0}^{\infty }} {\displaystyle \frac {
\mathrm{sin}(s\,\mathrm{arctan}(2\,t))\,e^{\pi \,t}}{(1 + 4\,t
^{2})^{ s/2 }\,\mathrm{cosh}(\pi \,t)}} \,dt \,. \hspace*{2cm} \sigma>1
\]
}
\label{NewJensen}
\end{equation}

Notice that (\ref{Sr2.3.41}) and (\ref{EtaC=half}) are valid for all $s$. Simple substitution of $c\rightarrow0$ in (\ref{Eta1WithC}) gives

\begin{equation}
\mapleinline{inert}{2d}{eta(s) = -sin(1/2*s*Pi)*Int(1/(sinh(Pi*t)*t^s),t = 0 .. infinity);}{%
\[
\eta (s)= - \mathrm{sin}({\displaystyle {s\,\pi /2 }} )\,
{\displaystyle \int _{0}^{\infty }} {\displaystyle \frac {t^{-s}}{
\mathrm{sinh}(\pi \,t)}} \,dt \hspace*{2cm} \sigma < 0
\]
}
\label{EtaC0}
\end{equation}

which, after taking (\ref{Eta=Zeta}) into account, leads to

\begin{equation}
\mapleinline{inert}{2d}{test := eta(s) = sin(1/2*s*Pi)/(2^(1-s)-1)*Int(t^{-s}/(sinh(Pi*t)),t = 0 .. infinity);}{%
\[
\zeta (s)={\displaystyle \frac {\mathrm{sin}(
{\displaystyle {s\,\pi \, /2 }} )}{(2^{1 - s} - 1)}} \,
{\displaystyle \int _{0}^{\infty }} {\displaystyle \frac {t^{-s}}{
\mathrm{sinh}(\pi \,t)\,}} \,dt, \hspace*{1cm} \sigma<0
\]
}
\label{ZetaC0}
\end{equation}

displaying the existence of the trivial zeros of $\zeta(s)$ at $s=-2\,N$, since the integral in (\ref{ZetaC0}) is clearly finite at each of these values of $s$. The result (\ref{ZetaC0}) is equivalent to the known result (\ref{Erd1.12.6}) after taking (\ref{functional}) into account. Other possibilities abound. For example, set $c=1/3$ in (\ref{Zeta1WithC}) to obtain
\begin{equation}
\zeta (s)=3^{s - 1}\,{\displaystyle \int _{0}^{\infty }} 
{\displaystyle \frac { (2\,
\mathrm{sin}(s\,\mathrm{arctan}(t))\,\mathrm{sinh}(
{\displaystyle 2\,\pi \,t\,/3} ) - \sqrt{3}\,\mathrm{cos}(
s\,\mathrm{arctan}(t)))}{(1 + t^{2})^{s/2} \, ( 2\,\mathrm{cosh}({\displaystyle 2
\,\pi \,t\,/3} ) + 1)}} \,dt \;. \hspace*{1cm} \sigma>1
\end{equation}

Further, a novel representation can be found by noting that 
\eqref{ContIntEta1} is invariant under the transformation of variables $t\rightarrow t+\rho$ together with the choice $c=\sigma$ along with the requirement that $0 \leq \sigma \leq 1$. This effectively allows one to choose $c=s$, yielding the following complex representation for $\eta(s)$, as well as the the corresponding representation for $\zeta(s)$ because of \eqref{Eta=Zeta}, valid in the critical strip $0\leq\sigma\leq 1$:

\begin{eqnarray}
{ \eta (s)={\displaystyle \frac {i}{2}} \,{\displaystyle \int _{s
 - i\infty }^{s + i\infty }} {\displaystyle \frac {t^{-s} }{
\mathrm{sin}(\pi \,t)}} \,dt }
 \label{ContIntEta0_with_c=s} \\ 
{\, \,\, ={\displaystyle \frac {1}{2}} \,{\displaystyle \int _{ - 
\infty }^{\infty }} {\displaystyle \frac {(s + i\,t)^{-s}}{
\mathrm{sin}(\pi \,(s + i\,t))}} \,dt \, \, .
\label{ContIntEta_with_c=s}
}
\end{eqnarray}

Because \eqref{ArcTanSin} and \eqref{ArcTanCos} are valid for $c \in \mathfrak{C}$, it follows that \eqref{Eta1WithC} with $c=s$ is also a valid representation for $\eta(s)$, as well as the the corresponding representation for $\zeta(s)$ with recourse to \eqref{Eta=Zeta}, again only valid in the critical strip $0 \leq \sigma \leq 1$. Another interesting variation is the choice $c=1/\rho$ in \eqref{Eta1WithC}, valid for $\rho>1$, again applicable to the corresponding representation for $\zeta(s)$ because of \eqref{Eta=Zeta}. Such a representation may find application in the analysis of $\zeta(\sigma+i\rho)$ in the asymptotic limit $\rho\rightarrow\infty$. Finally, the choice $c=\sigma,\,  0\leq\sigma\leq1$ leads to an interesting variation that is discussed in Section \ref{sec:crits}.

%% file: more_simple.tex
As a function of $c$, $\zeta(s)$ and $\eta(s)$ are constant, as are their integral representations over a range corresponding to $0<N-c<1$; however the integrals change discontinuously at each new value of $c$ that corresponds to a new value of $N$. Thus, the operator $\frac {\partial} {\partial c} $ acting on (\ref{ContInt}) and (\ref{ContIntEta1}) vanishes, except at non-negative integers $c=N$, when it becomes indefinite. Therefore, it is of interest to investigate this operation in the neighbourhood of $c=N/2$, as well as the limit $c \rightarrow N$. 

Applied to (\ref{ContIntEta1})\footnote{and equivalent to integration by parts}, the requirement that
\begin{equation} 
\partial {\eta(s)} / \partial {c} =0 \hspace*{2cm} c\approx N/2
\label{partial_eta}
\end{equation}

and identification of one of the resulting terms, yields

\begin{equation}
\mapleinline{inert}{2d}{eta(s+1) =
-1/2*1/s*Int((c+t*I)^(-s)/sin(Pi*(c+t*I))^2*cos(Pi*(c+t*I))*Pi,t =
-infinity .. infinity);}{%
\[
\eta (s+1 )= - {\displaystyle \frac {\pi}{2 \, s}} \,  \! 
{\displaystyle } \,{\displaystyle \int _{ - \infty }
^{\infty }} {\displaystyle \frac {(c + i\,t)^{- s}\,\mathrm{
cos}(\pi \,(c + i\,t)) }{\mathrm{sin}^{2}(\pi \,(c + i\,t)) }
} \,dt \!  ,
\]
}
\label{Eta(s+1)}
\end{equation}
a new result with improved convergence at infinity.  Similarly\footnote{It should be noted that Srivastava, Adamchik and Glasser (\cite{SG&A}) study an incomplete form of (\ref{Zetas+1}) in order to obtain representations for sums both involving $\zeta(n)$, and $\zeta(n)$ by itself, without commenting on the significance of the contour integral with infinite bounds (\ref{Zetas+1}).}
, operating with $\frac {\partial} {\partial c} $ on (\ref{ContInt})
yields

\begin{equation}
\zeta (s+1 )= {\displaystyle \frac {\pi}{2 \, s}} \,  \! 
{\displaystyle } \,{\displaystyle \int _{ - \infty }
^{\infty }} {\displaystyle \frac {(c + i\,t)^{- s}\,\mathrm{
} }{\mathrm{sin}^{2} (\pi \,(c + i\,t)) }
} \,dt \! \, .
\label{Zetas+1}
\end{equation}

Using (\ref{ArcTanSin}) and (\ref{ArcTanCos}) to reduce the integration range of (\ref{Eta(s+1)}), gives a fairly lengthy result for general values of $c$. If $c=1/2$, (\ref{Eta(s+1)}) becomes

\begin{equation}
\mapleinline{inert}{2d}{eta(s+1) =
Pi*2^(s-1)/s*Int((1+t^2)^(-1/2*s)*sinh(1/2*Pi*t)*sin(s*arctan(t))/cosh
(1/2*Pi*t)^2,t = 0 .. infinity);}{%
\[
\eta (s + 1)={\displaystyle \frac {2^{s - 1} \, \pi }{s}} \,
{\displaystyle \int _{0}^{\infty }} {\displaystyle \frac {(1 + t
^{2})^{- {s/2}}\,\mathrm{sinh}({\displaystyle 
\pi \,t/2} )\,\mathrm{sin}(s\,\mathrm{arctan}(t))}{\mathrm{cosh
}^{2}({\displaystyle {\pi \,t/2}} ) }} \,dt
\]
}
\label{Etas+1}
\end{equation}

which, by letting $s \rightarrow s-1$ can be rewritten as

\begin{equation}
\mapleinline{inert}{2d}{eta(s+1) =
2*(s-1)/s*eta(s)+Pi*2^(s-1)/s*Int(cos(s*arctan(t))*t*(1+t^2)^(-1/2*s)*
sinh(1/2*Pi*t)/cosh(1/2*Pi*t)^2,t = 0 .. infinity);}{%
\[
\eta (s + 1)={\displaystyle \frac {2\,(s - 1)\,\eta (s)}{s}}  + 
{\displaystyle \frac {2^{s - 1}\,\pi}{s}} \,{\displaystyle 
\int _{0}^{\infty }} {\displaystyle \frac {\mathrm{cos}(s\,
\mathrm{arctan}(t))\,(1 + t^{2})^{ - s/2}\,\mathrm{
sinh}({\displaystyle \pi \,t/2} )}{\mathrm{cosh}^{2} (
{\displaystyle \pi \,t/2} ) }} \, t \, dt  ,
\]
}
\label{Eta(s+1)_via_Eta(s)}
\end{equation}

and, with regard to (\ref{Eta=Zeta})
\begin{equation}
\mapleinline{inert}{2d}{Zeta(s+1) =
2*(s-1)/s*(1-2^(-s+1))*Zeta(s)/(1-2^(-s))+Pi*2^(s-1)/s/(1-2^(-s))*Int(
cos(s*arctan(t))*t*(1+t^2)^(-1/2*s)*sinh(1/2*Pi*t)/cosh(1/2*Pi*t)^2,t
= 0 .. infinity);}{%
\maplemultiline{
\zeta (s + 1)={\displaystyle \frac {2\,(s - 1)\,(1 - 2^{1 - s })\,\zeta (s)}{s\,(1 - 2^{ - s})}}  
\mbox{} \\ \hspace*{1cm} + {\displaystyle \frac {2^{s - 1}\,\pi}{s\,(1 - 2^{
 - s})}} \,{\displaystyle \int _{0}^{\infty }} {\displaystyle 
\frac {\mathrm{cos}(s\,\mathrm{arctan}(t))\,(1 + t^{2})^{-
s/2}\,\mathrm{sinh}({\displaystyle \pi \, t/2} 
)}{\mathrm{cosh}^{2} ({\displaystyle \pi \,t/2} ) }} \,t\,dt \,\, .
 }
}
\label{Zeta(s+1)_via_Zeta(s)}
\end{equation}

Similarly\footnote{ \eqref{Eta(s+1)_via_Eta(s)} and \eqref{Zeta(s+1)_via_Zeta(s)} are the first two entries in a family of recursive relationships that will be discussed elsewhere (in preparation).}
, in the case that $c=1/2$, (\ref{Zetas+1}) reduces to the well-known result (\cite{Lindelof}, Eq. 4,{\it{III}}(5)\footnote{attributed to Jensen \cite{Jensen}})

\begin{equation}
\mapleinline{inert}{2d}{Zeta(s+1) =
Pi/s*2^(s-1)*Int(cos(s*arctan(t))/((1+t^2)^(1/2*s))/cosh(1/2*Pi*t)^2,t
= 0 .. infinity);}{%
\[
\zeta (s + 1)={\displaystyle \frac {\pi \,2^{s - 1}}{s}} \,
{\displaystyle \int _{0}^{\infty }} {\displaystyle \frac {
\mathrm{cos}(s\,\mathrm{arctan}(t))}{(1 + t^{2})^{s/2}
\,\mathrm{cosh} ^{2} ({\displaystyle {\pi \,t}{/2} } )}} \,dt\;.
\]
}
\end{equation}

If $c \rightarrow 0$ in (\ref{Eta(s+1)})\footnote{For the case $c=0$ in (\ref{Zetas+1}) see the following section} , one finds, a finite result\footnote{also equivalent to integrating (\ref{EtaC0}) by parts}, 

\begin{equation}
\mapleinline{inert}{2d}{eta(s) =
Pi*sin(1/2*s*Pi)/(s-1)*Int(t^(-s+1)*cosh(Pi*t)/sinh(Pi*t)^2,t = 0 ..
infinity);}{%
\[
\eta (s)={\displaystyle \frac {\pi \,\mathrm{sin}({\displaystyle 
\frac {s\,\pi }{2}} )}{s - 1}} \,{\displaystyle \int _{0}^{\infty
 }} {\displaystyle \frac {t^{1 - s }\,\mathrm{cosh}(\pi \,t)
}{\mathrm{sinh}^{2} (\pi \,t) }} \,dt ,  \hspace*{2 cm} \sigma<0
\]
}
\label{EtaWithC=0ByParts}
\end{equation}

and, when $c=1$, one finds the unusual forms
\begin{equation}
\mapleinline{inert}{2d}{eta(s+1) =
1/2*I/s*Pi*x*Int((1+x*exp(theta*I)*I)^(-s)*cosh(Pi*x*exp(theta*I))*exp
(theta*I)/sinh(Pi*x*exp(theta*I))^2,theta = 0 ..
Pi)-x*Pi/s*Int((1+x^2*v^2)^(-1/2*s)*cos(s*arctan(x*v))*cosh(Pi*x*v)/si
nh(Pi*x*v)^2,v = 1 .. infinity);}{%
\maplemultiline{
\eta (s + 1)={\displaystyle \frac { 
\,i\,\pi \,r}{2\,s}} \,{\displaystyle \int _{0}^{\pi }} 
{\displaystyle \frac {(1 + r\,i\,e^{i\,\theta })^{( - s)}\,
\mathrm{cosh}(\pi \,r\,e^{i\,\theta })\,e^{i\,\theta }}{
\mathrm{sinh}^{2} (\pi \,r\,e^{i \,\theta }) }} \,d\theta  \\
\mbox{} -   \! {\displaystyle \frac {r\,\pi }{s}} \,
{\displaystyle \int _{1}^{\infty }} {\displaystyle \frac {(1 + r
^{2}\,t^{2})^{( - s/2)}\,\mathrm{cos}(s\,\mathrm{arctan}
(r\,t))\,\mathrm{cosh}(\pi \,r\,t)}{\mathrm{sinh}^{2}(\pi \,r\,t) 
}} \,dt \!  }
}
\label{Eta(s+1)_from_c=1}
\end{equation}

or, 
\begin{equation}
\mapleinline{inert}{2d}{eta(s+1) =
1/2*I/s*Pi*x*Int((1+x*exp(theta*I)*I)^(-s)*cosh(Pi*x*exp(theta*I))*exp
(theta*I)/sinh(Pi*x*exp(theta*I))^2,theta = 0 ..
Pi)-x*Pi/s*Int((1+x^2*v^2)^(-1/2*s)*cos(s*arctan(x*v))*cosh(Pi*x*v)/si
nh(Pi*x*v)^2,v = 1 .. infinity);}{%
\maplemultiline{
\eta (s + 1)=1+{\displaystyle \frac { 
\,i\,\pi \,r}{2\,s}} \,{\displaystyle \int _{0}^{\pi }} 
{\displaystyle \frac {(1 - r\,i\,e^{i\,\theta })^{( - s)}\,
\mathrm{cosh}(\pi \,r\,e^{i\,\theta })\,e^{i\,\theta }}{
\mathrm{sinh}^{2} (\pi \,r\,e^{i \,\theta }) }} \,d\theta  \\
\mbox{} -   \! {\displaystyle \frac {r\,\pi }{s}} \,
{\displaystyle \int _{1}^{\infty }} {\displaystyle \frac {(1 + r
^{2}\,t^{2})^{( - s/2)}\,\mathrm{cos}(s\,\mathrm{arctan}
(r\,t))\,\mathrm{cosh}(\pi \,r\,t)}{\mathrm{sinh}^{2}(\pi \,r\,t) 
}} \,dt \!  }
}
\label{Eta(s+1)_a_from_c=1}
\end{equation}

valid for any\footnote{particularly $r=1/\pi$ or $r=1/2$ } $0<r<1$ and all $s$, after deformation of the contour of integration in 
(\ref{Eta(s+1)}) either to the left or right of the singularity at $t=1$, respectively.

%% file: section_c=0.tex
With respect to (\ref{Zetas+1}), setting $c=0$ together with some simplification gives

\begin{equation}
{ \zeta(s) = -\pi^{s-1} \frac{\mathrm{sin} (\pi\,s/2)}{s-1} } \int _{0}^{\infty} \frac {t^{1-s}}{\mathrm{sinh}^{2}( t) } dt \hspace*{2cm} \sigma<0 \,
\label{Zeta_sinh^2}
\end{equation}

reproducing (\ref{GR527.1}), with $a=1$, after taking (\ref{functional}) into account.  This same case ($c=0$) naturally suggests that the integral in (\ref{ContInt0}) be broken into two parts

\begin{equation}
{\zeta(s) = \zeta_{0}(s)+\zeta_{1}(s)}
\label{ZetaParts}
\end{equation}

where

\begin{equation}
\zeta_{0}(s) =\frac{i}{2} \oint_{-i}^{i}\,t^{-s}\,\mathrm{cot}(\pi t)\, dt \,\,,
\label{Zeta_0}
\end{equation}

\begin{equation}
{\zeta_{1}(s) = \frac{i}{2} \int_{-i\,\infty}^{i}\,t^{-s}\,\mathrm{cot}(\pi t)\, dt
+\frac{i}{2} \int_{i}^{i\,\infty}\,t^{-s}\,\mathrm{cot}(\pi t)\, dt}
\label{ZetaC=0} 
\end{equation}

such that the contour in (\ref{Zeta_0}) bypasses the origin to the right. 
It is easy to show that the sum of the two integrals in (\ref{ZetaC=0}) reduces to

\begin{align}
{\zeta _{1}(s)} & = \mathrm{sin}({\displaystyle {\pi \,s /2}} )
\,{\displaystyle \int _{1}^{\infty }} t^{- s}\,\mathrm{coth}(
\pi \,t)\,dt  \hspace*{3cm} \sigma >1  \\
&={\displaystyle \frac {\mathrm{sin}({\displaystyle 
{\pi \,s/2}} )}{s - 1}}  +\, \mathrm{sin}({\displaystyle 
{\pi \,s/2}} )\,{\displaystyle \int _{1}^{\infty }} 
{\displaystyle \frac {t^{ - s}\,e^{( - \pi \,t)}}{\mathrm{sinh}
(\pi \,t)}} \,dt\,,
\label{zeta_1}
\end{align}

the latter valid for all $s$. One, of several commonly used possibilities of reducing (\ref{Zeta_0}) is to apply (\cite{Hansen} (54.3.2))

\begin{equation}
t\,\mathrm{cot}(\pi t)=-\frac{2}{\pi} \sum_{k=0}^{\infty}\,t^{2k}\,\zeta (2k) \hspace*{2cm} |t|<1
\label{Hansen_544.3.2}
\end{equation}  

giving

\begin{equation}
\zeta_{0}(s)=-\frac{2}{\pi }\,\mathrm{sin}(\frac{\pi s}{2}) \sum_{k=0}^{\infty}\,(-1)^{k}\,\frac{\zeta (2k)}{2k-s}\,\,,
\label{Sum1_s}
\end{equation}

convergent for all $s$. Such sums appear frequently in the literature (e.g. \cite{Crandall}, \cite{Sriv} Section (3.2)), and have been studied extensively when $s=-n$ (\cite{Sriv}, \cite{ChoiSrivAdam}). In the case $s=2n$, \eqref{Sum1_s} reduces to $\zeta_{0}(s)=\zeta(s)$ and \eqref{zeta_1} vanishes. \newline

Along with (\ref{Zeta_sinh^2}) the integral in (\ref{Zetas+1}) can be similarly split, using the convergent series representation 
\begin{equation}
\mapleinline{inert}{2d}{1/(sin^{2}(Pi*t)) = 2*Sum((2*k-1)*Zeta(2*k)*t^(2*k-2)/Pi^2,k = 0 ..
infinity);}{%
\[
{\displaystyle \frac {1}{\mathrm{sin}^{2} (\pi \,v) }} =2/{\pi^{2} }\,  \! {\displaystyle \sum _{k=0}^{\infty }} \,{\displaystyle 
{(2\,k - 1)\,\zeta (2\,k)\,v^{2\,k - 2}} } \,\,\,\,\,\,\,\,\, |v|<1
\]
}
\label{Hans50.6.10}
\end{equation}
 
where the Bernoulli numbers $B_{2k}$ appearing in the original (\cite{Hansen}, Eq.(50.6.10)) are related to $\zeta(2k)$ by

\begin{equation}
\mapleinline{inert}{2d}{Zeta(2*m) = -1/2*(-1)^m*B[2*m]*(2*Pi)^(2*m)/GAMMA(2*k+1);}{%
\[
\zeta (2\,k)= - {\displaystyle \frac {1}{2}} \,{\displaystyle 
\frac {(-1)^{k}\,  \,(2\,\pi )^{2\,k}\,{B_{2\,k}}}{\Gamma (2\,k + 1
)}} \,\,.
\]
}
\label{BernZ}
\end{equation}

Thus the representation
\begin{equation}
\mapleinline{inert}{2d}{Zeta(s) = -1/2*i*Pi/(s-1)*Int(t^(1-s)/sin(Pi*t)^2,t = -I ..
I)-Pi*sin(1/2*Pi*s)/(s-1)*Int(t^(1-s)/sinh(Pi*t)^2,t = 1 ..
infinity);}{%
\[
\zeta (s)= - {\displaystyle \frac {1}{2}} \,\! 
{\displaystyle \frac {i\,\pi }{(s - 1)}} \,{\displaystyle \int _{
 - i}^{i}} \,\, {\displaystyle \frac {t^{1 - s}}{\mathrm{sin} ^{2} (\pi \,
t)}} \,dt \!   -   \! {\displaystyle \frac {
\pi \,\mathrm{sin}({\displaystyle \pi \,s\,/2} )}{s - 1}} 
\,{\displaystyle \int _{1}^{\infty }} {\displaystyle \frac {t^{1
 - s}}{\mathrm{sinh}^{2} (\pi \,t) }} \,dt \!  
\]
}
\end{equation}

becomes (formally)
\begin{equation}
\mapleinline{inert}{2d}{Zeta(s) = sin(1/2*Pi*s)*(-2/Pi*Sum((-1)^k*(2*k-1)*Zeta(2*k)/(2*k-s),k
= 0 .. infinity)-Pi*Int(t^(1-s)/sinh(Pi*t)^2,t = 1 .. infinity));}{%
\[
\zeta (s)={\displaystyle - \frac {\mathrm{sin}({\displaystyle \pi \,s\,/2} )} {s-1}} \,
 \left(  \!  {\displaystyle \frac{2}{\pi} \,\,  
{\displaystyle \sum _{k=0}^{\infty }} \,
{\displaystyle \frac {(-1)^{k}\,(2\,k - 1)\,\zeta (2\,k)}{2\,k - s}}  \!  
}  + \pi \,{\displaystyle \int _{1}^{\infty }} {\displaystyle 
\frac {t^{1 - s}}{\mathrm{sinh}^{2} (\pi \,t) }} \,dt \! 
 \right) 
\]
}
\label{sinh^2TwoParts}
\end{equation}

after straightforward integration of (\ref{Hans50.6.10}) with $ \sigma < 0 $. In (\ref{sinh^2TwoParts}), the integral converges, and the formal sum conditionally converges or diverges according to several tests.  Howver, further analysis is still possible. Following the method of \cite{Sriv}, Section (3.3) applied to (\ref{Sum1_s}), together with \cite{Sriv}, Eq. 3.4(15), we have the following general result for $|v|<1$
\begin{eqnarray}
\mapleinline{inert}{2d}{Sum(Zeta(2*k)*t^(2*k)/(2*k-s),k = 0 .. infinity) =
1/(2*s)+Int(v^{-s}*Sum(Zeta(2*k)*v^(2*k-1),k = 1 .. infinity),v = 0 ..
t);}{%
\[
{\displaystyle \sum _{k=0}^{\infty }} \,{\displaystyle \frac {
\zeta (2\,k)\,v^{2\,k}}{2\,k - s}} ={\displaystyle \frac {1}{2
\,s}}  + {\displaystyle \int _{0}^{v}} t^{- s}\,
{\displaystyle \sum _{k=1}^{\infty }} \,\zeta (2\,k)\,t^{2\,k - 
1}\,dt
\label{PsiDiff0}
\]
} \\
\mapleinline{inert}{2d}{Sum(Zeta(2*k)*t^(2*k)/(2*k-s),k = 0 .. infinity) =
1/(2*s)+1/2*t^s*Int(t^(-s)*(Psi(1+t)-Psi(1-t)),t = 0 .. infinity);}{%
\[
={\displaystyle \frac {1}{2
\,s}}  + {\displaystyle \frac {1}{2}} \,v^{s}\,{\displaystyle 
\int _{0}^{v}} t^{- s}\,(\Psi (1 + t) - \Psi (1 - t))\,dt \hspace*{1cm} \sigma<2 \, 
\]
}
\label{PsiDiff}
\end{eqnarray}
so that (formally), after premultiplying \eqref{PsiDiff} with $1/v$, then operating with $v^2 \frac {\partial}{\partial v}$,  
 
\begin{eqnarray}
\mapleresult
\mapleinline{inert}{2d}{Sum(Zeta(2*k)*(2*k-1)*t^(2*k)/(2*k-s),k = 0 .. infinity) =
-1/(2*s)+1/2*t^s*Int(v^(1-s)*(Psi(1,1+v)+Psi(1,1-v)),v = 0 .. t);}{%
\[
{\displaystyle \sum _{k=0}^{\infty }} \,{\displaystyle \frac {
\zeta (2\,k)\,(2\,k - 1)\,v^{2\,k}}{2\,k - s}} = - 
{\displaystyle \frac {1}{2\,s}}  + {\displaystyle \frac {1}{2}} 
\,v^{s}\,{\displaystyle \int _{0}^{v}} t^{1 - s}\,(\Psi^{\,\prime} ( \,1 + t) + \Psi^{\, \prime} ( \,1 - t))\,dt
\]
} \\
\mapleinline{inert}{2d}{Eq := Sum(Zeta(2*k)*(2*k-1)*t^(2*k)/(2*k-s),k = 0 .. infinity) =
-1/(2*s)+1/2*t^s*Int(v^(1-s)*(Pi^2*(1+cot(Pi*v)^2)-1/(v^2)),v = 0 ..
t);}{%
\[
= - {\displaystyle \frac {1}{2\,s}}  + {\displaystyle 
\frac {1}{2}} \,v^{s}\,{\displaystyle \int _{0}^{v}} t^{1 - s}
\,(\pi ^{2}\,(1 + \mathrm{cot}^{2} (\pi \,t) ) - {\displaystyle 
\frac {1}{t^{2}}} )\,dt . \hspace*{1.5 cm} \sigma<2 \, 
\]
}
\label{IntCot^2}
\end{eqnarray} 

Notice that the right-hand sides of (\ref{PsiDiff}) provides an analytic continuation of the left-hand side for $|v|>1$ and the right-hand side of (\ref{IntCot^2}) gives a regularization of the left-hand side. Putting together (\ref{sinh^2TwoParts}) and (\ref{IntCot^2}),  incorporating the first term of (\ref{IntCot^2}) into the second term of (\ref{sinh^2TwoParts}), setting $v=i$ and constraining $\sigma>0$ eventually yields a new integral representation\footnote{compare with \cite{TyHo}, Eq. (2.8).} 

\begin{equation}
\mapleinline{inert}{2d}{Zeta(s) =
-Pi^(s-1)*sin(1/2*Pi*s)/(s-1)*Int(v^(1-s)*(1/(sinh(v)^2)-1/(v^2)),v =
0 .. infinity);}{%
\[
\zeta (s)= -  \! {\displaystyle \frac {\pi ^{s - 1}\,
\mathrm{sin}({\displaystyle \frac {\pi \,s}{2}} )}{s - 1}} \,
{\displaystyle \int _{0}^{\infty }} t^{1 - s}\,({\displaystyle 
\frac {1}{\mathrm{sinh}^{2}\, t }}  - {\displaystyle \frac {1}{t^{2
}}} )\,dt \!   
\]
}
\label{Int_sinh^2a}
\end{equation}

valid for $0<\sigma<2$, a non-overlapping analytic continuation of (\ref{Zeta_sinh^2}) that spans the critical strip. Additionally, using (\ref{functional}) followed by replacing $s\rightarrow 1-s$ gives a second representation 
\begin{equation}
\mapleinline{inert}{2d}{Zeta(s) = 2^(s-1)/GAMMA(1+s)*Int(t^s*(1/(sinh(t)^2)-1/(t^2)),t = 0 ..
infinity);}{%
\[
\zeta (s)={\displaystyle \frac {2^{s - 1}}{\Gamma (1 + s)}} \,
{\displaystyle \int _{0}^{\infty }} t^{s}\,({\displaystyle 
\frac {1}{\mathrm{sinh} ^{2}\, t}}  - {\displaystyle \frac {1}{t^{2
}}} )\,dt
\]
}
\label{Int_sinh^2b}
\end{equation}

valid for $-1<\sigma<1$, a non-overlapping analytic continuation of (\ref{GR527.1}). Over the critical strip $0<\sigma<1$ representations (\ref{Int_sinh^2a}) and (\ref{Int_sinh^2b}) are simultaneously true and augment the rather short list of other representations that share this property (\cite{Bruijn}, as cited in \cite{Bateman}\footnote{unnumbered equations immediately following 1.12(1), also reproduced in \cite{SrivChoi2}, Eqs. 2.3(44)-(46).}). See Section \ref{sec:crits} for further application of this result.

%% file: c_geq_1.tex
Interesting results can be obtained by translating the contour in (\ref{ContInt0}) and (\ref{ContIntEta1}) and adding the residues so-omitted. From (\ref{ZetaSum}) and (\ref{EtaSum}) define the partial sums

\begin{equation}
\zeta_{N}(s) \equiv \sum_{k=1}^{N} {k^{-s}} 
\label{ZetaSumN}
\end{equation}

\begin{equation}
\eta_{N}(s) \equiv -\sum_{k=1}^{N} {(-1)^k\,k^{-s}} 
\label{EtaSumN},
\end{equation}

let $c \rightarrow c_{N}, N = 1,2,3,... $ with 
\begin{equation}
N\leq c_{N}< N+1
\label{cLimits}
\end{equation}

and adjust\footnote{If $c_{N}=N$ an extra half-residue is removed and the contour deformed to pass to the right of the point $t=N$.} all the above results accordingly by redefining $c\rightarrow c_{N}$. This gives, for example, the generalization of (\ref{ContInt0})

\begin{equation}
{\zeta(s) = \zeta_{N}(s) -\frac{1} {2N^s} \delta _{N,c_{N}} -\frac{1}{2} \int_{c_{N}-i\infty}^{c_{N}+i\infty}\,t^{-s}\,cot(\pi t)\, dt} \hspace*{.3cm} \;,\;\sigma>1
\label{ContInt0N} 
\end{equation} 

and (\ref{ContIntEta1})

\begin{equation}
{ \eta (s)=\eta_{N}(s) +\frac {(-1)^{N} } {2N^s} \delta_{N,c_{N}}  + {\displaystyle \frac {1}{2}} \,{\displaystyle \int _{c_{N}
 - i\infty }^{c_{N} + i\infty }} {\displaystyle \frac {t^{-s} }{
\mathrm{sin}(\pi \,t)}} \,dt \, \, }
\label{ContIntEta}
\end{equation}

so (\ref{Zeta1WithC}) and (\ref{Eta1WithC}) with $c \rightarrow c_{N} $ become integral representations of the remainder of the partial sums $\zeta(s)-\zeta_{N}(s)$ and $\eta(s)-\eta_{N}(s)$, the former corresponding to a special case of the Hurwitz zeta function $\zeta(s,N)$.\newline

In the case that $c_{N}=N=1$, one finds

\begin{align}
\zeta(s) & = {\displaystyle \frac {1}{2}}  + {\displaystyle \int _{0
}^{\infty }} {\displaystyle \frac {\mathrm{cosh}(\pi \,t)\,
\mathrm{sin}(s\,\mathrm{arctan}(t))}{(1 + t^{2})^{s/2}
\,\mathrm{sinh}(\pi \,t)}} \,dt \hspace*{3cm} \sigma>1
\label{ZetaC1}
\\ 
&= {\displaystyle \frac {s + 1}{2\,(
s - 1)}}  + {\displaystyle \int _{0}^{\infty }} {\displaystyle 
\frac {e^{ - \pi \,t}\,\mathrm{sin}(s\,\mathrm{arctan}(t))}{(1
 + t^{2})^{s/2}\,\mathrm{sinh}(\pi \,t)}} \,dt 
\label{ZetaC1a}
\end{align}

and

\begin{equation}
\mapleinline{inert}{2d}{eta(s) = 1/2+Int(sin(s*arctan(t))/sinh(Pi*t)/((1+t^2)^(1/2*s)),t = 0
.. infinity);}{%
\[
\eta (s)={\displaystyle \frac {1}{2}}  + {\displaystyle \int _{0}
^{\infty }} {\displaystyle \frac {\mathrm{sin}(s\,\mathrm{arctan}
(t))}{\mathrm{sinh}(\pi \,t)\,(1 + t^{2})^{s/2}}} \,dt
\]
}
\label{EtaC1}
\end{equation}

with the latter two results valid for all $s$. Neither of \eqref{ZetaC1} nor \eqref{EtaC1} appear in references; \eqref{ZetaC1a} corresponds to a limiting case of Hermite's representation of the Hurwitz zeta function (\cite{Bateman}, Eq. 1.10(7)). In analogy to the developments leading to (\ref{NewJensen}), equations (\ref{ZetaC1}), (\ref{EtaC1}) and (\ref{Eta=Zeta}), together with the following modified form of (\cite{G&R}, Eq.(9.51.5)) 

\begin{equation}
\mapleinline{inert}{2d}{zeta(s) =
2^(s-1)*s/(2^s-1)/(s-1)+2^(s-1)/(2^s-1)*Int(sin(s*arctan(t))*(1+exp(-P
i*t))/((1+t^2)^(1/2*s))/sinh(Pi*t),t = 0 .. infinity);}{%
\[
\zeta (s)={\displaystyle \frac {2^{s - 1}\,s}{(2^{s} - 1)\,(s
 - 1)}}  + {\displaystyle \frac {2^{s - 1}}{2^{s} - 1}} \,
{\displaystyle \int _{0}^{\infty }} {\displaystyle \frac {
\mathrm{sin}(s\,\mathrm{arctan}(t))\,(1 + e^{ - \pi \,t})}{(1
 + t^{2})^{ s/2 }\,\mathrm{sinh}(\pi \,t)}} \,dt
\]
}
\end{equation}

can be used to obtain the possibly new result

\begin{equation}
\mapleinline{inert}{2d}{zeta(s) =
1/2*(s-3)/(s-1)+Int(sin(s*arctan(v))*exp(Pi*v)/((1+v^2)^(1/2*s))/sinh(
Pi*v),v = 0 .. infinity);}{%
\[
\zeta (s)={\displaystyle \frac {s - 3}{2\,(s - 1)}}  + 
{\displaystyle \int _{0}^{\infty }} {\displaystyle \frac {
\mathrm{sin}(s\,\mathrm{arctan}(t))\,e^{\pi \,t}}{(1 + t^{2})^{
s/2}\,\mathrm{sinh}(\pi \,t)}} \,dt \; .\hspace*{2cm} \; \sigma>1
\]
}
\end{equation}

For $c=1/2+N$ corresponding to general partial sums, two possibilities are apparent. Based on (\ref{ArcTanSin}), (\ref{ArcTanCos}), (\ref{Zeta1WithC}) and (\ref{Eta1WithC})\footnote{for comparison, see also \cite{Olver}, Section 8.3, and \cite{Edwards}, Section (6.4) and Chapter 7.} we find 

\begin{align}
\zeta (s) -  {\displaystyle \sum _{k=1}^{N}} \,
{\displaystyle \frac {1}{k^{s}}}& ={\displaystyle 
\frac {1}{2}} \,({\displaystyle 1/2}  + N)^{1 - s }
\,{\displaystyle \int _{0}^{\infty }} {\displaystyle \frac {
\mathrm{sinh}(\pi \,t\,(1 + 2\,N))\,(1 + t^{2})^{- s/2}\,\mathrm{sin}(s\,\mathrm{arctan}(t))}{\mathrm{cosh} ^{2} (
{\displaystyle {\pi \,t\,(1/2 + N)}})}} \,dt \hspace*{1cm} \sigma>1
\label{partial sum1} 
\\
& ={\displaystyle \frac {({\displaystyle 1/2}  + N)^{1 - s
}}{s - 1 }}  - 2\,({\displaystyle 1/2}  + N)^{(1 - s)
}\,{\displaystyle \int _{0}^{\infty }} {\displaystyle \frac {(1
 + t^{2})^{- s/2}\,\mathrm{sin}(s\,\mathrm{arctan}(t)
)}{1 + e^{\pi \,t\,(1 + 2\,N)}}} \,dt
\label{partial sum2}
\end{align}
and
\begin{equation}\mapleinline{inert}{2d}{eta(s)+Sum((-1)^k/(k^s),k = 1 .. N) =
(-1)^N*(1/2+N)^(-s+1)*Int((1+t^2)^(-1/2*s)*cos(s*arctan(t))/cosh(Pi*(1
/2+N)*t),t = 0 .. infinity);}{%
\[
\eta (s) +  {\displaystyle \sum _{k=1}^{N}} \,
{\displaystyle \frac {(-1)^{k}}{k^{s}}}  =(-1)^{N}\,(
{\displaystyle 1/2}  + N)^{1 - s}\,{\displaystyle 
\int _{0}^{\infty }} {\displaystyle \frac {(1 + t^{2})^{ - 
s/2}\,\mathrm{cos}(s\,\mathrm{arctan}(t))}{\mathrm{cosh
}(\pi \,t\,({\displaystyle 1/2}  + N))}} \,dt
\]
}
\end{equation}

which together extend (\ref{JensenZ}) and (\ref{EtaC=half}) respectively. Application of (\ref{Eta(s+1)}) and (\ref{Zetas+1}) with (\ref{ArcTanSin}) and (\ref{ArcTanCos}), also with $c=1/2+N$, gives the remainders

\begin{equation}
\mapleinline{inert}{2d}{Zeta(s+1)-Sum(1/(k^(s+1)),k = 1 .. N) =
Pi*(1/2+N)^(-s+1)/s*Int((1+t^2)^(-1/2*s)*cos(s*arctan(t))/cosh(Pi*(1/2
+N)*t)^2,t = 0 .. infinity);}{%
\[
\zeta (s + 1) -  {\displaystyle \sum _{k=1}^{N}} \,
{\displaystyle \frac {1}{k^{s + 1}}} =
{\displaystyle \frac {\pi \,({\displaystyle 1/2}  + N)^{
1 - s }}{s}} \,{\displaystyle \int _{0}^{\infty }} 
{\displaystyle \frac {(1 + t^{2})^{ - s/2}\,\mathrm{
cos}(s\,\mathrm{arctan}(t))}{\mathrm{cosh}^{2} (\pi \,t\,({\displaystyle 
1/2 }  + N)) }} \,dt
\]
}
\label{Zeta with c=N/2}
\end{equation}
and
\begin{equation}
\mapleinline{inert}{2d}{eta(s+1)+Sum((-1)^k/(k^(s+1)),k = 1 .. N) =
Pi*(-1)^N*(1/2+N)^(-s+1)/s*Int(sinh(Pi*(1/2+N)*t)*(1+t^2)^(-1/2*s)*sin
(s*arctan(t))/cosh(Pi*(1/2+N)*t)^2,t = 0 .. infinity);}{%
\maplemultiline{
\eta (s + 1) +  {\displaystyle \sum _{k=1}^{N}} \,
{\displaystyle \frac {(-1)^{k}}{k^{s + 1}}} = \\
{\displaystyle \frac {\pi \,(-1)^{N}\,({ 1/2}  + N)^{1 - s }}{s}} \,{\displaystyle \int _{0}^{\infty }} 
{\displaystyle \frac {\mathrm{sinh}(\pi \,t\,({ 1/2}  + N))\,(1 + t^{2})^{- s/2}\,\mathrm{sin}(s
\,\mathrm{arctan}(t))}{\mathrm{cosh} ^{2} (\pi \,t\,({ 
1/2 }  + N))}} \,dt }
}
\label{Eta with c=N/2}
\end{equation}
valid for all $s$. For the case $c=N$ applying (\ref{ContInt0N}) and (\ref{ContIntEta}) gives the remainders
\begin{align}
\zeta (s) -  {\displaystyle \sum _{k=1}^{N}} \,
{\displaystyle \frac {1}{k^{s}}}& = - {\displaystyle 
\frac {1}{2\,N^{s}}}  + N^{1 - s }\,{\displaystyle \int _{0}
^{\infty }} {\displaystyle \frac {\mathrm{cosh}(\pi \,N\,t)\,(1
 + t^{2})^{- s/2}\,\mathrm{sin}(s\,\mathrm{arctan}(t)
)}{\mathrm{sinh}(\pi \,N\,t)}} \,dt \hspace*{1cm} \sigma>1\\
&= - {\displaystyle \frac {1}{2\,N^{s}}}  + {\displaystyle \frac {N
^{1 - s}}{ s - 1 }}  + 2\,N^{1 - s}\,{\displaystyle \int _{0
}^{\infty }} {\displaystyle \frac {(1 + t^{2})^{ - s/2
}\,\mathrm{sin}(s\,\mathrm{arctan}(t))}{e^{2\,N\,\pi \,t} - 1}
} \,dt
\label{Zeta_with_c=N}
\end{align}

and 

\begin{equation}
\eta (s) +  {\displaystyle \sum _{k=1}^{N}} \,
{\displaystyle \frac {(-1)^{k}}{k^{s}}} =
{\displaystyle \frac {(-1)^{N}}{2\,N^{s}}}  + (-1)^{1 + N}\,N^{
1 - s }\,{\displaystyle \int _{0}^{\infty }} {\displaystyle 
\frac {(1 + t^{2})^{- s/2}\,\mathrm{sin}(s\,\mathrm{
arctan}(t))}{\mathrm{sinh}(\pi \,N\,t)}} \,dt 
\end{equation}
which extend (\ref{ZetaC1}) and (\ref{EtaC1}), the latter valid for all $s$. It should be noted that \eqref{Zeta_with_c=N}, valid for $s \in \mathfrak{C}$ generalizes well-established results (e.g. \cite{G&R}, Eq. 3.411(12)) that are only valid in the limiting case $s=n$\footnote{Although \cite{G&R}, Eq. 3.411(21) gives (with a missing minus sign) an integral representation for finite sums such as those appearing in \eqref{Zeta_with_c=N}, specified to be valid only for integer exponents ($k^n$), that particular result in \cite{G&R} also appears to apply after the generalization $n \rightarrow s$.}. Omitting the integral term from \eqref{partial sum2} gives an upper bound for the remainder; similarly, omitting the integral term in \eqref{Zeta_with_c=N} gives a lower bound for the remainder, provided that $s \in \mathfrak{R}$ in both cases. In addition, many of the above, in particular \eqref{partial sum2} and \eqref{Zeta_with_c=N} appear to be new and should be compared to classical results such as \cite{Ivic}, Theorem 1.8, page 23 and pages 99-100.

%% file: by_parts.tex
By simple manipulation of the integral representations given, it is possible to obtain new series representations. Starting from the integral representation (\ref{ZetaC0}), integrate by parts, giving

\begin{equation}
\mapleinline{inert}{2d}{Zeta(s) =
sin(1/2*Pi*s)/(2^(1-s)-1)*Pi^2/(-1/2*s+1)*Int(t^(-s+2)*hypergeom([-1/2
*s+1],[3/2, 2-1/2*s],1/4*t^2*Pi^2)/sinh(Pi*t)^3*cosh(Pi*t),t = 0 ..
infinity);}{%
\maplemultiline{
\zeta (s)={\displaystyle \frac {\mathrm{sin}({\displaystyle 
\pi \,s/\,2} )\,\pi ^{2}}{(2^{(1 - s)} - 1)\,( - 
{\displaystyle s\,/2}  + 1)}}  \\
\times {\displaystyle \int _{0}^{\infty }} {\displaystyle \frac {t^{2- 
s}\,\mathrm{_{1}F_{2}}( - {\displaystyle \frac {s}{2}}  + 1
; \,{\displaystyle \frac {3}{2}} , \,2 - {\displaystyle \frac {
s}{2}} ; \,{\displaystyle \frac {t^{2}\,\pi ^{2}}{4}} )\,
\mathrm{cosh}(\pi \,t)}{\mathrm{sinh}^{3} (\pi \,t) }} \,dt \hspace*{3cm} \sigma<0 \,\, ,
} 
\label{Int1F2}
}
\end{equation}

where the hypergeometric function arises because of \cite{Luke}, Eq.6.2.11(6). Write the hypergeometric function as a (uniformly convergent) series, thereby permitting the sum and integral to be interchanged, apply a second integration by parts, and find

\begin{equation}
\mapleinline{inert}{2d}{Zeta(s) =
sin(1/2*Pi*s)*Pi*Sum(Pi^(2*k)/GAMMA(2*k+2)*Int(1/sinh(Pi*t)^2*t^(-s+1+
2*k),t = 0 .. infinity),k = 0 .. infinity)/(2^(1-s)-1);}{%
\[
\zeta (s)={\displaystyle \left( \frac {\pi}{2^{1 - s} - 1} \right )} \,{\displaystyle \mathrm{sin}({\displaystyle 
\pi \,s/2} )\, \! {\displaystyle \sum _{k=
0}^{\infty }} \, {\displaystyle \frac {\pi ^{2\,k}}{
\Gamma (2\,k + 2)}} \,{\displaystyle \int _{0}^{\infty }} 
{\displaystyle \frac {t^{1 - s  + 2\,k}}{\mathrm{sinh}^{2} (\pi \,
t) }} \,dt  } 
\hspace*{3cm} \sigma<0 \, \, .
\]
\label{SumInt}
}
\end{equation}

Identify the integral in (\ref{SumInt}) using (\ref{GR527.1}) giving
\begin{equation}\mapleinline{inert}{2d}{Zeta(s) =
sin(1/2*Pi*s)*2^s*Pi^(s-1)*Sum(2^(-2*k)*Zeta(-s+1+2*k)/GAMMA(2*k+2)*GA
MMA(-s+2+2*k),k = 0 .. infinity)/(2^(1-s)-1);}{%
\[
\zeta (s)={\displaystyle \frac {2^{s}\,\pi ^{s - 1}\, \mathrm{sin}({\displaystyle 
\pi\,s/2 } )\,}{2^{1 - s} - 1} \, 
{\displaystyle \sum _{k=0}^{\infty }} \,{\displaystyle \frac {2^{
- 2\,k}\,\zeta ( - s + 1 + 2\,k)\,\Gamma ( - s + 2 + 2\,k)}{
\Gamma (2\,k + 2)}}  }   
\]
\label{ZGamSum}
}
\end{equation}
valid for all values of $\sigma$ by the ratio test and the principle of analytic continuation. This result, which appears to be new, can be rewritten by extracting the $k=0$ term, and applying (\ref{functional}) to yield

\begin{equation}\mapleinline{inert}{2d}{Zeta(s) =
sin(1/2*Pi*s)*2^s*Pi^(s-1)*Sum(2^(-2*k)*Zeta(-s+1+2*k)/GAMMA(2*k+2)*GA
MMA(-s+2+2*k),k = 1 .. infinity)/(2^(1-s)-2+s);}{%
\[
\zeta (s)={\displaystyle \frac {2^{s}\,\pi ^{s - 1}\,\mathrm{sin}({\displaystyle \pi\,s/2 } )\,  \! 
} {2^{1 - s} - 2 + s}} {\displaystyle \sum _{k=1}^{\infty }} \,{\displaystyle \frac {2^{
- 2\,k}\,\zeta ( - s + 1 + 2\,k)\,\Gamma ( - s + 2 + 2\,k)}{
\Gamma (2\,k + 2)}}   
\],
\label{SumInt2}
}
\end{equation}

a representation, valid for all values of $s$, that is reminiscent of the following similar sum recently obtained by Tyagi and Holm (\cite{TyHo} Eq.(3.5))

\begin{equation}
\zeta (s) = \frac{\pi^{s-1}}{1-2^{1-s}} \mathrm{sin}( \pi\,s/2)\sum _{k=1} ^{\infty} \frac{(2-2^{s-2k})\,\zeta(-s+1+2k)\Gamma(-s+1+2k) }{\Gamma(2k+2)} \hspace*{2cm}  \sigma>{0}\,  \footnote{by Gauss' test}
\label{TyHo}
\end{equation}

Eq. (\ref{SumInt2}) can also be rewritten in the form
\begin{equation}
\zeta(s)=\frac{1}{\Gamma(s)(2^{s}-1-s)}\sum _{k=1} ^{\infty} \frac{2^{-2k} \zeta(s+2k)\Gamma(1+s+2k)}{\Gamma(2k+2)}
\label{SumInt2a}
\end{equation}
by replacing $s\rightarrow 1-s$ and applying (\ref{functional}).\newline

With reference to (\ref{ContInt}), express the numerator $\it{cosine}$ factor as a difference of exponentials each of which in turn is written as a convergent power series in $\pi\,(c+it)$, recognize that the resulting series can be interchanged with the integral and each term identified using (\ref{ContIntEta1}) to find 

\begin{equation}
\mapleinline{inert}{2d}{zeta(s) =
-Pi^(1/2)*Sum((-1/4*Pi^2)^k/GAMMA(k+1/2)/GAMMA(k+1)*eta(s-2*k),k = 0
.. infinity);}{%
\[
\zeta (s)= - \sqrt{\pi }\, \left(  \! {\displaystyle \, \sum _{k=0}
^{\infty }} \,{\displaystyle \frac {( - {\displaystyle {\pi
 ^{2}/4}} )^{k}\,\eta (s - 2\,k)}{\Gamma (k + {1/2} )\,\Gamma (k + 1)}}  \!  \right)\,.  \hspace*{2cm} \sigma >1 
\]
\label{Zsum1}
}
\end{equation}

A more interesting variation of (\ref{Zsum1}) can be obtained by applying (\ref{Eta=Zeta}) and explicitly identifying the $k=0$ term of the sum giving 

\begin{equation}
\mapleinline{inert}{2d}{zeta(s) =
-1/2*Pi^(1/2)*Sum((-1/4*Pi^2)^k/GAMMA(k+1/2)/GAMMA(k+1)*(1-2^(1-s+2*k)
)*Zeta(s-2*k),k = 1 .. infinity)/(1-2^(-s));}{%
\[
\zeta (s)= - {\displaystyle \frac {1}{2}} \,{\displaystyle 
\frac { \! {\displaystyle \sum _{k=1}^{
\infty }} \,{\displaystyle \frac {( - {\displaystyle {\pi 
^{2}}} )^{k}\,(1 - 2^{1 - s + 2\,k}) }{
\Gamma (2k + 1)}\zeta (s - 2\,k) }  \! 
 }{1 - 2^{- s}}}, 
\]
\label{Zsum1a}
}
\end{equation}

a convergent representation if $\sigma >1$ (by Gauss' test). An immediate consequence is the identification

\begin{equation}
\zeta(2) = - \frac {\pi^{2}} {3} \zeta(0) .
\end{equation}

because only the $k=1$ term in the sum contributes when $s=2$. Additionally, when $s=2n$, (\ref{Zsum1a}) terminates at $n$ terms, becoming a recursion for $\zeta(2n)$ in terms of $\zeta(2n-2k), k=1, \dots n$ - see section \ref{sec:s=n} for further discussion. If the functional equation (\ref{functional}) is used to extend its region of convergence, (\ref{Zsum1a}) becomes

\begin{equation}
\mapleinline{inert}{2d}{zeta(s) =
1/2*Pi^s*Sum(Zeta(1-s+2*k)*(2^(s-2*k-1)-1)/GAMMA(2*k+1)/GAMMA(s-2*k),k
= 1 .. infinity)/(-1+2^(-s))/cos(1/2*Pi*s);}{%
\[
\zeta (s)={\displaystyle \frac {1}{2}} \,{\displaystyle \frac {
\pi ^{s}\, {\displaystyle \sum _{k=1}^{\infty }} \,
{\displaystyle \frac {\zeta (1 - s + 2\,k)\,(2^{s - 2\,k - 1}
 - 1)}{\Gamma (2\,k + 1)\,\Gamma (s - 2\,k)}}  }{( 2^{- s}-1)\,\mathrm{cos}({\displaystyle {\pi \,s/\,2}} 
)}} \hspace*{1cm} \sigma>1
\]
}
\end{equation}

which augments similar series appearing in \cite{Sriv}, Sections 4.2 and 4.3, and corresponds to a combination of entries in \cite{Hansen} (Eqs. (54.6.3) and (54.6.4)).

%% file: split2.tex
\subsubsection{The lower range}\label{sect::lower}

After applying (\ref{functional}), the integral in (\ref{Zeta_sinh^2}) is frequently (e.g. \cite{SrivChoi2}) and conveniently split in two, defining

\begin{equation}
\zeta_{-}(s) = \int_{0}^{\;\omega} \frac {t^s}{\mathrm{sinh}^2(t)} dt \qquad  \sigma>1
\label{Zeta-}
\end{equation}

and

\begin{equation}
\zeta_{+}(s) = \int_{\omega}^{\infty} \frac {t^s}{\mathrm{sinh}^2(t)} dt
\label{Zeta+}
\end{equation}

with

\begin{equation}
\zeta(s)\frac{\Gamma(s+1)}{2^{s-1}} = \zeta_{-}(s) +\zeta_{+}(s)
\label{ZetaSplit}  
\end{equation}

for arbitrary $\omega$. Similar to the derivation of (\ref{sinh^2TwoParts}), after application of (\ref{functional}) we have

\begin{equation}
\mapleinline{inert}{2d}{Zeta[x] =
-2*Pi^(1/2)*Sum(bernoulli(2*k)/GAMMA(k-1/2)/GAMMA(k+1)/(2*k-1+s),k = 0
.. infinity);}{%
\[
{\zeta _{-}(s)}= - 2\,\sqrt{\pi }\;  \! {\displaystyle \sum 
_{k=0}^{\infty }} \,{\displaystyle \frac {\mathrm{B}_{2\,k
} \;\omega^{2k+s-1} }{\Gamma (k - {\displaystyle 1/\,2} )\,\Gamma (k + 1)\,(2
\,k - 1 + s)}}  \!  
\]
}
\end{equation}

which can be rewritten using (\ref{BernZ}) as

\begin{equation}
\mapleinline{inert}{2d}{Zeta[x] =
4*Sum((-1/(Pi^2))^k*GAMMA(k+1/2)*Zeta(2*k)/GAMMA(k-1/2)/(2*k-1+s),k =
0 .. infinity);}{%
\[
{\zeta _{-}(s)}=4\, \omega^{s-1} \! {\displaystyle \sum _{k=0}^{\infty }
} \,{\displaystyle \frac {{( - {\displaystyle \frac {\omega^{2}}{\pi ^{2}}
} )}^{k}\,\Gamma (k + {\displaystyle \frac {1}{2}} )\,\zeta (2\,k)
}{\Gamma (k - {\displaystyle 1/\,2} )\,(2\,k - 1 + s)}}\,. \hspace*{1cm} \omega\leq\pi 
 \!   
\]
}
\label{zeta_{-}1}
\end{equation}

Substitute (\ref{ZetaSum}) in (\ref{zeta_{-}1}), interchange the order of summation, and identify the resulting series to obtain

\begin{equation}
\mapleinline{inert}{2d}{Jbb5 := 2*s/Pi^2*Sum(hypergeom([1,
1/2+1/2*s],[3/2+1/2*s],-1/(Pi^2*m^2))/m^2,m = 1 ..
infinity)/(s+1)-coth(1)+s/(s-1);}{%
\[
{\zeta _{-}(s)}= {\displaystyle \frac {2\,s\,\omega^{s+1}} {\pi ^{2}\,(s + 1)}  \! 
{\displaystyle \sum _{m=1}^{\infty }} \,{\displaystyle \frac {1} {m^{2}} \;\,
\mathrm{_{2}{ \displaystyle F}_{1}}(1, \,{\displaystyle \frac {1}{2}}  + 
{\displaystyle \frac {s}{2}} ; \,{\displaystyle \frac {3}{2}} 
 + {\displaystyle \frac {s}{2}} ; \, - {\displaystyle \frac {\omega^ 2  }{
\pi ^{2}\,m^{2}}} )}   \!   } 
 - \omega^{s}\mathrm{coth}(\omega) + {\displaystyle \frac {s\,\omega^{s-1}}{s - 1}}\, , 
\]
}
\label{zeta_{-}2}
\end{equation}

a convergent representation that is valid for all $s$.\newline

Other representations are easily obtained by applying any of the standard hypergeometric linear transforms as in \eqref{Ftransform} - see below. One such leads to

\begin{equation}
\mapleinline{inert}{2d}{Zeta(s) =
s*GAMMA(1/2+1/2*s)*Sum(GAMMA(k+1)*Sum((1/(Pi^2*m^2+1))^(k+1),m = 1 ..
infinity)/GAMMA(k+3/2+1/2*s),k = 0 .. infinity)-coth(1)+s/(s-1);}{%
\[
\zeta_{-}(s)=s\,\omega^{s-1}\,\Gamma ({\displaystyle \frac {1}{2}}  + 
{\displaystyle \frac {s}{2}} )\,  \! {\displaystyle \sum 
_{k=0}^{\infty }} \,{\displaystyle \frac {\Gamma (k + 1) } {\Gamma (k + {\frac {3}{2}}  + 
{ \frac {s}{2}} )} \,
 \! {\displaystyle \sum _{m=1}^{\infty }} \,(
{\displaystyle \frac {1}{\pi ^{2}\,m^{2}/\omega^{2} + 1}} )^{k + 1} \! 
}  \!   - \omega^{s}\mathrm{coth}(\omega)
 + {\displaystyle \frac {s\,\omega^{s-1}}{s - 1}} 
\]
}
\label{zeta_{-}3}
\end{equation}

a rapidly converging representation for small values of $\omega$, valid for all $s$. The inner sum of (\ref{zeta_{-}3}) can be written in several ways. By expanding the denominator as a series in $1/(t^2m^2)$ with $r<2p-1$ and interchanging the two series we find the general form 

\begin{equation}
\mapleinline{inert}{2d}{Sum(m^r/((1+t^2*m^2)^p),m = 1 .. infinity) =
Sum(GAMMA(p+m)*(-1)^m*t^(-2*m-2*p)*Zeta(2*m+2*p-r)/GAMMA(p)/GAMMA(m+1)
,m = 0 .. infinity);}{%
\[
{\kappa(r,p,t) \equiv \displaystyle \sum _{m=1}^{\infty }} \,{\displaystyle \frac {m^{
r}}{(1 + t^{2}\,m^{2})^{p}}} ={\displaystyle \sum _{m=0}^{\infty 
}} \,{\displaystyle \frac {\Gamma (p + m)\,(-1)^{m}\,t^{- 2\,m
 - 2\,p}\,\zeta (2\,m + 2\,p - r)}{\Gamma (p)\,\Gamma (m + 1)}} 
\]
}
\label{xform1}
\end{equation}

where the right hand side continues the left if $|t|>1$ and $r-2\,p\geq -1$. Alternatively, \cite{Hansen}, Eq. (6.1.32) identifies

\begin{equation}
\sum_{m=0}^{\infty}\,\frac{1}{(m^2\,t^2+y^2)}=\frac{1}{2y^2}+\frac{\pi}{2\,t\,y}\mathrm{coth}(\frac{\pi y}{t})
\label{Hansen6.1.32}
\end{equation}

so for $p=k+1$ and $r=0$ in (\ref{xform1}), each term of the inner series (\ref{zeta_{-}3}) can be obtained analytically by taking the $k^{th}$ derivative of the right-hand side with respect to the variable $y^2$ at $y=1$, leading to a rational polynomial in coth$(\omega)$ and sinh$(\omega)$. The first few terms are given below for two choices of $\omega$. \newline

\begin{tabular}[t] {|l|ccc|}
\multicolumn{4} {c} {The first few values of $\kappa(0,k+1,\pi^2/\omega^2)$ in (\ref{xform1}) } \vspace*{0.2cm} \\\hline
\it{k}& analytic&\ $\omega =1$ &\ $ \omega=\pi$ \\\hline &&& \\ 
0 & $  - 1/2   + 
\,\omega \,\mathrm{coth}(\omega)/2 $ & 0.156518&1.076674\\&&&\\
1 &$ - {\displaystyle \frac {1}{2}}  + {\displaystyle \frac {1}{4}} 
\,{\displaystyle \frac {\omega ^{2} }{\mathrm{sinh}^{2}\omega  }
}  + {\displaystyle \frac {1}{4}} \,\omega \,\mathrm{coth}(\omega
 )
$&$0.927424\times10^{-2}$&0.306837\\ & & & \\ 
2 &$ \frac {3\,\omega }{16} \,\mathrm{coth}^{3}\omega -  \frac {1}{2} \, \mathrm{coth}^{2} \omega  + 
\, \frac {\omega \,
\mathrm{coth}\,\omega\,( - 3 + 2\,\omega ^{2})}{16\,\mathrm{sinh}^{2}
\omega }$  &&\\  & $+ \frac {8 + 3\,\omega ^{2}}{16\,\mathrm{sinh}^{2}\omega} 
$&$0.795491\times10^{-3}$&0.134296\\&&&
\\\hline
\end{tabular}\vspace*{0.4cm}

Application of (\ref{xform1}) to (\ref{zeta_{-}3}) leads to a less rapidly converging representation that can be expressed in terms of more familiar functions

\begin{equation}
\mapleinline{inert}{2d}{Zeta(s) =
s*Sum(Sum((-1)^m*GAMMA(m+k+1)*Zeta(2*m+2*k+2)/GAMMA(m+1)/(Pi^(2*m)),m
= 0 .. infinity)/(Pi^(2*k+2))/GAMMA(k+3/2+1/2*s),k = 0 ..
infinity)*GAMMA(1/2+1/2*s)-coth(1)+s/(s-1);}{%
\maplemultiline{
\zeta_{-}(s)=s\, \omega^{s-1} \!\,\Gamma ({\displaystyle \frac {1
}{2}}  + {\displaystyle \frac {s}{2}} )  {\displaystyle \sum _{k=0}^{\infty }} \,
{\displaystyle {{\displaystyle \sum _{m=0}^{\infty }} \,
{\displaystyle \frac {(-1)^{m}\,\Gamma (m + k + 1)\,\zeta (2\,m
 + 2\,k + 2)}{ \Gamma (k + { \frac {3}{2}}  + { 
\frac {s}{2}} )    \Gamma (m + 1)\, }} } 
{ }} {\displaystyle {(\frac {\omega^{2}}{\pi^{2}})} ^ {(m +k+1)}} \\ \hspace*{1cm} - \omega^{s}\mathrm{coth}(\omega)
 + { \displaystyle \frac {s\,\omega^{s-1}}{s - 1}}\;. }}
\label{zeta_{-}5} 
\end{equation}


A second transform of the intermediate hypergeometric equation yields another useful result
\begin{equation}
\mapleinline{inert}{2d}{Zeta[x] =
2*s/GAMMA(1/2+1/2*s)*Sum(Pi^(-2*k)*Sum((-1)^m*GAMMA(m+k+1/2+1/2*s)*Zet
a(2*m+2*k+2)/GAMMA(m+1)/(Pi^(2*m)),m = 0 ..
infinity)/(2*k+1+s)/GAMMA(k+1),k = 0 ..
infinity)/Pi^2-coth(1)+s/(s-1);}{%
\maplemultiline{
{\zeta _{-}(s)}={\displaystyle \frac {2\,s\,\omega^{s-1}} {\Gamma ({\displaystyle 
1/2 }  + {\displaystyle s/2} ) }} \: 
{\displaystyle \sum _{k=0}^{\infty }} \,{\displaystyle {{( {\frac{\omega^{2}}{\pi^{2}}})}^{( k+1)}\, \! {\displaystyle \sum _{m=0}^{\infty }} {(\frac{\omega^2}{\pi^2})}^{m}
\,{\displaystyle \frac {(-1)^{m}\,\Gamma (m + k + {\displaystyle 
\frac {1}{2}}  + {\displaystyle \frac {s}{2}} )\,\zeta (2\,m + 2
\,k + 2)}{\Gamma (m + 1)\,\Gamma (k + 1)\,(2\,k + 
1 + s)\,  }}  \!  }{}} \\ 
\hspace*{1cm} - \omega^{s}\mathrm{coth}(\omega)
 + { \displaystyle \frac {s\,\omega^{s-1}}{s - 1}} \; .} }
\label{zeta_{-}6} 
\end{equation}

Both (\ref{zeta_{-}5}) and (\ref{zeta_{-}6} ) can be re-ordered along a diagonal of the double sum\footnote{$ \sum_{k=0}^{\infty} \sum_{m=0}^{\infty}f(m,k) =\sum_{k=0}^{\infty} \sum_{m=0}^{ k }f(m,k-m)  $}, allowing one of the resulting (terminating) series to be evaluated in closed form. In either case, both eventually give

\begin{equation}
\mapleinline{inert}{2d}{Zeta(s) = -2*s*Sum((-1)^k*Zeta(2*k)/(Pi^(2*k))/(2*k-1+s),k = 0 ..
infinity)-coth(1);}{%
\[
\zeta_{-}(s)= - 2\,s\,\omega^{s-1} \! {\displaystyle \sum _{k=0}^{\infty
 }} \,{\displaystyle {(\frac{-\omega^{2}}{\pi^2})}^k \frac {\zeta (2\,k)}{
\,(2\,k - 1 + s)}}  \!   - \omega^{s}\,\mathrm{coth}(\omega)\,,
\]
\label{zeta_{-}7} 
}
\end{equation}

reducing to (\ref{Sum1_s}) in the case $\omega = \pi$ taking (\ref{functional}) into consideration. Furthermore, by applying (\ref{ZetaSum}) and interchanging the resulting series, a further (more slowly converging) representation is found

\begin{equation}
\zeta_{-}(s)=-2 \, \frac{\omega^{s-1}}{s-1}\sum_{m=1}^{\infty}
\mathrm{_{2}{ \displaystyle F}_{1}}
(1,\frac{s}{2}-\frac{1}{2}\, ; \frac{s}{2}+\frac{1}{2}; \frac{-\omega^2}{\pi^2 \,m^2}) - \omega^{s}\,\mathrm{coth}(\omega)\, ,
\label{zeta_{-}8}
\end{equation}

which differs from (\ref{zeta_{-}2})\footnote{ and could have also been obtained by recourse to contiguity relations.}. Series of the form (\ref{zeta_{-}7}) have been studied extensively (\cite{Sriv}, \cite{SrivChoi2}) when $s=n$ for special values of $\omega$. The representations given above will find application in Section \ref{sec:s=n}. \newline 

\subsubsection{The upper range}

Rewrite (\ref{Zeta+}) in exponential form and expand the factor $(1-\mathrm{e}^{-2t})^{-2}$ in a convergent power series, interchange the sum and integral and recognize that the resulting integral is a generalized exponential integral $\mathrm{E}_s(z)$ usually defined in terms of incomplete Gamma functions (\cite{Milg2}, \cite{OLBC}, Section 8.19 ) as follows
\begin{equation}
\mapleinline{inert}{2d}{Ei(s, z) = z^(s-1)*GAMMA(1-s, z)}{\[\displaystyle {\mathrm E_{ \, s}} \left( z \right) ={z}^{s-1}\Gamma  \left( 1-s,z \right) \]}
\end{equation}

giving a rapidly converging series representation

\begin{equation}
\zeta_{+}(s) =4\,\omega^{s+1}\,\sum_{k=1}^{\infty} k \, \mathrm E_{-s}(2\omega\,k\,).
\label{E_sSum}
\end{equation}

A useful variant of (\ref{E_sSum}) is obtained by manipulating the well-known recurrence relation (\cite{OLBC}, Eq. 8.19.12)  for $\mathrm{E}_s(z)$ to obtain

\begin{equation}
\mathrm{E}_{-s}(z) = \frac {1}{z^2} \left( \mathrm{e}^{-z}(s+z) \,+\,s(s-1)\,\mathrm{E}_{2-s}(z) \right)\,.
\label{recurE}
\end{equation}

Apply (\ref{recurE}) repeatedly $N$ times to (\ref{E_sSum}) and find

\begin{equation}
\mapleinline{inert}{2d}{Zeta(s) =
-omega^(-1+s)*s*ln(-1+exp(2*omega))+2*omega^s*s-2*omega^s/(1-exp(2*ome
ga))+2*omega^s*GAMMA(s+1)*Sum(Ly[k](exp(-2*omega))/GAMMA(s-k+1)/((2*om
ega)^k),k = 2 ..
N+1)+2*omega^(s-N-1)*GAMMA(s+1)*Sum(E[N+2-s](2*k*omega)/(k^(N+1)),k =
1 .. infinity)/GAMMA(s-N-1)/(2^(N+1));}{%
\maplemultiline{
\zeta_{+}(s)= - s\, \omega ^{ - 1 + s }\,\mathrm{log}( e^{2\,
\omega }- 1 ) + 2\,s\,\omega ^{s}\, - {\displaystyle \frac {2\,\omega 
^{s}}{1 - e^{2\,\omega }}}  
\mbox{} + 2\,\omega ^{s}\,\Gamma (s + 1)\,  \! 
{\displaystyle \sum _{k=2}^{N + 1}} \,{\displaystyle \frac {{
\mathrm{Li}_{k}}(e^{ - 2\,\omega })}{\Gamma (s - k + 1)\,(2\,
\omega )^{k}}}  \!   \\
\hspace*{1cm} \mbox{} + {\displaystyle \frac {2\,\omega ^{s - N - 1}\,\Gamma 
(s + 1)\, }{\Gamma (s - N - 1)\,2^{N + 1}}}  \! {\displaystyle \sum _{k=1}^{\infty }} \,
{\displaystyle \frac {{\mathrm{E}_{N + 2 - s}}(2\,k\,\omega )}{k^{N + 1}
}}  \!   }
}
\label{EiSum_with_N}
\end{equation}

in terms of logarithms and polylogarithms ($\mathrm{Li}_{k}(z)$), a result that diverges numerically as $N$ increases, but which also, usefully, terminates if $s=n, N \geq n\,\geq 0$. This is discussed in more detail in Section \ref{sec:s=n}.\newline

Another useful representation emerges from (\ref{E_sSum}) by writing the function $\mathrm{E}_{s}(..)$ in terms of confluent hypergeometric functions. From \cite{Milg2}, Eq.(2.6b) and \cite{Lebedev}, Eq.(9.10.3) we find 

\begin{equation}
\mapleinline{inert}{2d}{4*omega^(s+1)*sum(k*Ei(-s,2*k*omega),k = 1 .. infinity) =
2^(1-s)*GAMMA(s+1)*Zeta(s)-4*omega^{s+1}*Sum(k*exp(-2*k*omega)*hyperge
om([1],[2+s],2*k*omega),k = 1 .. infinity)/(s+1);}{%
\maplemultiline{
\zeta_{+}(s)\,=\,4\,\omega ^{s + 1}\,{\displaystyle \sum _{k=1}^{\infty }} \,k
\,\mathrm{E}_{-s}(2\,k\,\omega )= \\
2^{1 - s}\,\Gamma (s + 1)\,\zeta (s) - {\displaystyle \frac {4
\,\omega ^{s + 1}}{s + 1} \,{\displaystyle \sum _{k=1}^{\infty }} \,k\,
e^{ - 2\,k\,\omega }\,\mathrm{_{1}F_{1}}(1 ; \,2 + s; \,2\,k
\,\omega ) } , }
}
\label{Esum_as_1F1}
\end{equation}

where the series on the right hand side converges absolutely if $\sigma >1$ independent of $\omega$. Comparison of (\ref{E_sSum}), (\ref{ZetaSplit}) and (\ref{Esum_as_1F1}) identifies

\begin{equation}
\zeta_{-}(s)= {\displaystyle \frac {4
\,\omega ^{s + 1}}{s + 1} \,{\displaystyle \sum _{k=1}^{\infty }} \,k\,
e^{ - 2\,k\,\omega }\,\mathrm{_{1}F_{1}}(1 ; \,2 + s; \,2\,k
\,\omega ) }  
\end{equation}

which can be interpreted as a transformation of sums, or an analytic continuation, for the various guises of $\zeta_{-}(s)$ given in Section(\ref{sect::lower}), if both, or either, of the two sides converge for some choice of the variables $s$ and/or $\omega$, respectively. The fact that a term containing $\zeta(s)$ naturally appears in the expression (\ref{Esum_as_1F1}) for $\zeta_{+}(s)$ suggests that there will be a severe cancellation of digits between the two terms $\zeta_{-}(s)$ and $\zeta_{+}(s)$ if (\ref{ZetaSplit}) is used to attempt a numerical evaluation of $\zeta(s)$, as has been reported elsewhere (\cite{Dunster}). Essentially, (\ref{ZetaSplit}) is numerically useless (with common choices of arithmetic precision) for any choices of $\rho>10$, even if $\omega$ is carefully chosen in an attempt to balance the cancellation of digits. However, with 20 digits of precision, using 20 terms in the series, it is possible to find the first zero of $\zeta(1/2+i\rho)$ correct to 10 digits, deteriorating rapidly thereafter as $\rho$ increases. For any choice of $\omega$ , the convergence properties of the two sums representing $\zeta_{-}(s)$ and $\zeta_{+}(s)$ anti-correlate, as can be seen from the Table in the previous Section. This can also be seen by rewriting (\ref{Esum_as_1F1}) to yield the combined representation

\begin{equation}
\mapleinline{inert}{2d}{2^(1-s)*GAMMA(s+1)*Zeta(s) =
4*omega^(s+1)*Sum(k*(E_{-s}(2*k*omega)+hypergeom([s+1],[2+s],-2*k*omega
)/(s+1)),k = 1 .. infinity);}{%
\[
2^{1 - s}\,\Gamma (s + 1)\,\zeta (s)=4\,\omega ^{s + 1}\,
 \left(  \! {\displaystyle \sum _{k=1}^{\infty }} \,k\,(\mathrm{
E_{-s}}( 2\,k\,\omega ) + {\displaystyle \frac {\mathrm{
_{1}F_{1}}(s + 1; \,2 + s; \, - 2\,k\,\omega )}{s + 1}} ) \! 
 \right) \,. 
\]
}
\label{Esum_combined}
\end{equation}
In (\ref{Esum_combined}), the first inner term ($\mathrm{E}_{-s}(2\,k\,\omega)$) vanishes like $\mathrm{exp}(-2\,k\,\omega)/(2\,k\,\omega)$ whereas the second term $(_{1}F_{1}/(s+1))$ vanishes like $1/k^{(1+s)} - \mathrm{exp}(-2\,k\,\omega)/(2\,k\,\omega)$ so for large values of $k$ cancellation of digits is bound to occur\footnote{In fact, the coefficients  of both asymptotic series prefaced by the exponential terms cancel exactly term by term.}.\newline

It should be noted that \eqref{Esum_as_1F1} is representative of a family of related transformations. Although setting $s=0$ in \eqref{Esum_as_1F1} does not lead to a new result, operating first with $\frac{\partial}{\partial s}$ followed by setting $s=0$ leads to a transformation of a related sum
\begin{equation}
\mapleinline{inert}{2d}{Sum(E(1,2*k*omega),k = 1 .. infinity) =
-omega*arctanh((-omega^2)^(1/2)/Pi)/Pi/(-omega^2)^(1/2)+Sum(-omega*arc
tanh((-omega^2)^(1/2)/Pi/k)/Pi/k/(-omega^2)^(1/2),k = 2 ..
infinity)+1/2*ln(omega/Pi)+1/(2*omega)+1/2*gamma;}{%
\maplemultiline{
{\displaystyle \sum _{k=1}^{\infty }} \,\mathrm{E_1}(\,2\,k\,
\omega )= 
 - {\displaystyle \frac {\omega \,\mathrm{arctanh}(
{\displaystyle \frac {\sqrt{ - \omega ^{2}}}{\pi }} )}{\pi \,
\sqrt{ - \omega ^{2}}}}  -  {\displaystyle \sum _{k=2}
^{\infty }} {\displaystyle \frac {\omega \,
\mathrm{arctanh}({\displaystyle \frac {\sqrt{ - \omega ^{2}}}{\pi
 \,k}} )}{\pi \,k\,\sqrt{ - \omega ^{2}}}}  \!  
 + {\displaystyle \frac {1}{2}} \,\mathrm{log}(
{\displaystyle \frac {\omega }{\pi }} ) + {\displaystyle \frac {1
}{2\,\omega }}  + {\displaystyle \frac {\gamma }{2}}  }
}
\label{E1sum}
\end{equation}

where the first term corresponding to $k=1$ has been isolated, since it contains the singularity that occurs due to divergence of the left-hand sum if $\omega=i\,\pi$. Similarly, operating with $\frac{\partial}{\partial s}$ on \eqref{Esum_as_1F1} at $s=1$, leads to the identification

\begin{equation}
{\displaystyle \sum _{k=1}^{\infty }} \,\mathrm{log}(1 + 
{\displaystyle \frac {1}{\pi ^{2}\,k^{2}}} )= - 1 - \mathrm{log}(2
) + \mathrm{log}( e^{2} - 1 )
\label{logsum}
\end{equation}

(which can alternatively be obtained by writing the infinite sum as a product). A second application of $\frac{\partial}{\partial s}$ on \eqref{Esum_as_1F1}, again at $s=1$ leads to the following transformation of sums

\begin{equation}
\mapleinline{inert}{2d}{Sum(1/2*Zeta(2*k)*(-1)^k*Pi^(-2*k)/k^2,k = 1 .. infinity) =
1/12*Pi^2-1/2*gamma^2-ln(2)-gamma(1)+1/2*ln(2)^2+1/2-1/2*Sum(ln((Pi^2*
k^2+1)/Pi^2/k^2),k = 1 ..
infinity)-ln(-1/2*(-2+2*exp(2))^(1/2)/(exp(-2)-1))-Sum(Ei(1,2*k)/k,k =
1 .. infinity);}{%
\maplemultiline{
{\displaystyle \frac {1
}{2} \sum _{k=1}^{\infty }} \,{\displaystyle } \,{\displaystyle \frac {\zeta (2\,k)\,(-1)^{k}\,\pi ^{- 2
\,k}}{k^{2}}} ={\displaystyle \frac {\pi ^{2}}{12}}  - 
{\displaystyle \frac {\gamma ^{2}}{2}}  - \mathrm{log}(2) - \gamma
 (1) + {\displaystyle \frac {1}{2}} \,\mathrm{log^{2}}(2) + 
{\displaystyle \frac {1}{2}}  \\
\mbox{} - {\displaystyle \frac {1}{2}} \, 
{\displaystyle \sum _{k=1}^{\infty }} \,\mathrm{log}(
{\displaystyle \frac {\pi ^{2}\,k^{2} + 1}{\pi ^{2}\,k^{2}}} )
  - \mathrm{log}( - {\displaystyle \frac {1}{2}} \,
{\displaystyle \frac {\sqrt{ - 2 + 2\,e^{2}}}{e^{-2} - 1}} ) - 
 \! {\displaystyle \sum _{k=1}^{\infty }} \,
{\displaystyle \frac {\mathrm{E_1}(2\,k)}{k}}  \,\,.  
 }
}
\label{D2s=1}
\end{equation}

This transformation can be further reduced by applying an integral representation given in \cite{OLBC}, Section 6.2, leading to 

\begin{equation}
{\displaystyle \sum _{k=1}^{\infty }} \,{\displaystyle \frac {
\mathrm{E_1}(2\,k)}{k}} = - {\displaystyle \int _{1}^{\infty 
}} {\displaystyle \frac {\mathrm{log}(1 - e^{- 2\,t})}{t}} \,dt \,.
\label{E1(k)/k}
\end{equation}

Correcting a misprinted expression given in (\cite{Adam&Sriv}, Eq. (4.8), \cite{Sriv} and \cite{SrivChoi2}, Eq. 3.3(43)) which should read
\begin{equation}
\sum_{k=1}^{\infty} {\displaystyle \frac{t^{k}} {k^2}}\, \zeta(2\,k) = \int_{0}^{t}\, \mathrm{log}(\pi \, \sqrt{v} \: \mathrm{csc}( \pi \, \sqrt{v})\;)\,\frac{dv}{v},
\end{equation}
followed by a change of variables, eventually yields the representation
\begin{equation}
\mapleinline{inert}{2d}{Sum(Zeta(2*k)*(-1)^k*Pi^(-2*k)/k^2,k = 1 .. infinity) =
-2*Int((ln(1/2*exp(t)-1/2*exp(-t))-ln(t))/t,t = 0 .. 1);}{%
\[
{\displaystyle \sum _{k=1}^{\infty }} \,{\displaystyle \frac {
\zeta (2\,k)\,(-1)^{k}\,\pi ^{( - 2\,k)}}{k^{2}}} = - 2\,
{\displaystyle \int _{0}^{1}} {\displaystyle \frac {1} {t} ( {\mathrm{log}(
{\displaystyle \frac {1}{2}} \,e^{t} - {\displaystyle \frac {1}{2
}} \,e^{- t}) - \mathrm{log}(t)} } ) \,dt \,.
\]
}
\label{SumZ(2k)}
\end{equation}

Substituting \eqref{logsum}, \eqref{E1(k)/k}) and \eqref{SumZ(2k)} into \eqref{D2s=1} followed by some simplification, finally gives
\begin{equation}
\mapleinline{inert}{2d}{Int(ln(1/2*(1-exp(-2/v))*(exp(2*v)-1)/v)/v,v = 0 .. 1) =
-1/12*Pi^2+1/2*gamma^2+gamma(1)-1/2*ln(2)^2+2;}{%
\[
{\displaystyle \int _{0}^{1}} {\displaystyle {\mathrm{log}
 \left(  \! {\displaystyle} \,{\displaystyle \frac {
(1 - e^{- 2/v})\,(e^{2\,v} - 1)}{2\,v}}  \!  \right) 
} } \,{\displaystyle \frac {dv} {v}} = - {\displaystyle \frac {\pi ^{2}}{12}}  + 
{\displaystyle \frac {\gamma ^{2}}{2}}  + \gamma (1) - 
{\displaystyle \frac {1}{2}} \,\mathrm{log}^{2}(2) + 2\,,
\]
}
\end{equation}
a result that appears to be new.

\subsection{Another Interpretation}

These expressions for the split representation can be interpreted differently. For the moment, let $\sigma<-3$, and insert the contour integral representation \cite{Milg2}, Eq. (2.6a) for $\mathrm{E}_{-s}(z)$ into the converging series representation (\ref{E_sSum}), noting that contributions from infinity vanish, giving

\begin{equation}
\mapleinline{inert}{2d}{Zeta(s) = -Sum(-1/2*I/Pi*Int(GAMMA(-t)*(2*k*omega)^t/(s+1+t),t),k = 1
.. infinity);}{%
\[
\zeta_{+} (s)= - 4\,{\omega}^{s+1}\, {\displaystyle \sum _{k=1}^{\infty }} \,
 {\displaystyle  {{\displaystyle \frac {1}{2 \, \pi i}} 
}} \,{\displaystyle \int_{c-i\infty}^{c+i\infty} } {\displaystyle \frac {\Gamma 
( - t)\,(2\,k\,\omega )^{t}}{s + 1 + t}} \,dt  \! 
\]
}
\label{E_Contour1}
\end{equation}

where the contour is defined by the real parameter $c<-2$. With this caveat, the series in (\ref{E_Contour1}) converges, allowing the integration and summation to be interchanged, leading to the identification

\begin{equation}
\mapleinline{inert}{2d}{Zeta(s) =
2*I*omega^(s+1)/Pi*Int(GAMMA(-t)*(2*omega)^t*Zeta(-1-t)/(s+1+t),t =
c-infinity*I .. c+infinity*I);}{%
\[
\zeta_{+} (s)=- {\displaystyle \frac {4\,\omega ^{s + 1}}{2\,\pi\,i }} 
\,{\displaystyle \int _{c - i\,\infty } ^{c + i\,\infty }} 
{\displaystyle \frac {\Gamma ( - t)\,(2\,\omega )^{t}\,\zeta ( - 
1 - t)}{s + 1 + t}} \,dt\,.
\]
}
\label{E_Contour2}
\end{equation}

The  contour may now be deformed to enclose the singularity at $t=-s-1$, translated to the right to pick up residues of $\Gamma(-t)$ at $t=n, n=0,1\dots$ and the residue of $\zeta(-1-t)$ at $t=-2$. After evaluation, each of the residues so obtained cancels one of the terms in the series expression (\ref{zeta_{-}1}) for $\zeta_{-}(s)$ leaving (after a change of integration variables)

\begin{equation}
\mapleinline{inert}{2d}{GAMMA(s+1)*Zeta(s)/(2^(s-1)) =
1/2*I*omega^(s+1)/Pi*Int(GAMMA(1+t)*Zeta(t)/(s-t)/(2^(-1+t))/(omega^(1
+t)),t);}{%
\[
{\displaystyle \frac {\Gamma (s + 1)\,\zeta (s)}{2^{s - 1}\,\omega^{s+1}}} =
- {\displaystyle \frac {{\displaystyle 1}}{2\,\pi\,i }} \,{\displaystyle \oint } {\displaystyle \frac {
\Gamma (t + 1)\,\zeta (t)}{(s - t)\,2^{ t - 1 }\,\omega ^{t
 + 1}}} \,dt\, ,
\]
}
\end{equation} 

the standard Cauchy representation of the left-hand side, valid for all $s$, where the contour encloses the simple pole at $t=s$.

%% file: by_parts_2.tex
Integration by parts applied to previous results yields other representations. 

From (\ref{Erd1.12.4}), we obtain\footnote{known to Maple, but only in the form of a Mellin transform}

\begin{equation}
\zeta(s) = - \frac {1}{\Gamma(s-1)} \int_{0}^{\infty}\, t^{s-2}\, \mathrm{log}(1-\mathrm{exp}(-t))\, dt \hspace*{2cm} \sigma>1\,.
\end{equation}.\newline

and from (\ref{ZetaC1}) we find

\begin{equation}
\mapleinline{inert}{2d}{Zeta(s) =
1/2-1/Pi*s*Int(ln(sinh(Pi*t))*cos(arctan(t)*(s+1))/((1+t^2)^(1/2+1/2*s
)),t = 0 .. infinity);}{%
\[
\zeta (s)={\displaystyle \frac {1}{2}}  -   \! 
{\displaystyle \frac {s}{\pi }} \,{\displaystyle \int _{0}^{
\infty }} {\displaystyle \frac {\mathrm{log}(\mathrm{sinh}(\pi \,t
))\,\mathrm{cos}(\mathrm{(s+1)\,arctan}(t))}{(1 + t^{2})^{(1+s)/2}}} \,dt \! \hspace*{1 cm} \sigma>1 
\]
\label{Z320a1}
}
\end{equation}
as well as the following intriguing result\footnote{discovered by the Maple computer code}
\begin{equation}
\mapleinline{inert}{2d}{Zeta(s) =
1/2+1/2*Pi^(1/2)*sin(1/2*Pi*s)*GAMMA(1/2*s-1/2)/GAMMA(1/2*s)+Pi*Int(t/
sinh(Pi*t)^2*sin(s*arctan(t))*_{2}F_{1}([1/2, 1/2*s],[3/2],-t^2),t = 0
..
infinity)-s*Int(t*cosh(Pi*t)*cos(s*arctan(t))/sinh(Pi*t)/(1+t^2)*hyper
geom([1/2, 1/2*s],[3/2],-t^2),t = 0 .. infinity);}{%
\maplemultiline{
\zeta (s)={\displaystyle \frac {1}{2}}  + {\displaystyle \frac {1
}{2}} \,{\displaystyle \frac {\sqrt{\pi }\,\mathrm{sin}(
{\displaystyle \frac {\pi \,s}{2}} )\,\Gamma ({\displaystyle 
\frac {s}{2}}  - {\displaystyle \frac {1}{2}} )}{\Gamma (
{\displaystyle \frac {s}{2}} )}}  + \pi \,{\displaystyle \int _{0
}^{\infty }} {\displaystyle \frac {t\,\mathrm{sin}(s\,\mathrm{
arctan}(t))\, }{\mathrm{sinh}^{2} (\pi \,t) }} \mathrm{_{2}F_{1}}({\displaystyle 1/2} , 
\,{\displaystyle s/2} ; \,{\displaystyle 3/2
} ; \, - t^{2})  \,dt \\
\mbox{} - s\,{\displaystyle \int _{0}^{\infty }} {\displaystyle 
\frac {t\,\mathrm{cosh}(\pi \,t)\,\mathrm{cos}(s\,\mathrm{arctan}
(t))\,  } {\mathrm{sinh}(\pi \,t)\,(1 + t^{2})}} \mathrm{_{2}F_{1}}({\displaystyle 1/2} , \,
{\displaystyle s/2} ; \,{\displaystyle 3/2} 
; \, - t^{2}) \,dt \,. \hspace*{1 cm} \sigma>0 }
\label{Z320a2}
}
\end{equation}

This result has the interesting property that it is the only integral representation of which I am aware in which the independent variable $s$ does not appear as an exponent. Many variations of (\ref{Z320a2}) can be obtained through the use of any of the transformations of the hypergeometric function, for example

\begin{equation}
\mapleinline{inert}{2d}{hypergeom([1/2, 1/2*s],[3/2],-t^2) = 1/(1+t^2)^(1/2)*hypergeom([1/2,
3/2-1/2*s],[3/2],t^2/(1+t^2));}{%
\[
\mathrm{_{2}F_{1}}({\displaystyle 1/2 } , \,
{\displaystyle s/2} ; \,{\displaystyle 3/2} 
; \, - t^{2})={\displaystyle {({1 + t^{2}})^{-s/2}\,\, } { \mathrm{_{2}F_{1}}(
{\displaystyle 1} , \,{\displaystyle s/2} 
; \,{\displaystyle 3/2 }
; \,{\displaystyle {t^{2}}/(1 + t^{2}) } )}, 
} 
\]
}
\label{Ftransform}
\end{equation}
which choice happens to restore the variable $s$ to its usual role as an exponent.

%% file: int_reps.tex
From (\ref{SumInt2}), identify the product $\zeta(...)\Gamma(...)$ in the form of an integral representation from one of the primary definitions of $\zeta(s)$ (\cite{G&R}, Eq. (9.511)) and interchange the series and sum. The sums can now be explicitly evaluated, leaving an integral representation valid over an important range of $s$

\begin{equation}
\mapleinline{inert}{2d}{Zeta(s) =
Pi^(-1+s)*4^s*sin(1/2*Pi*s)/(2-2^(1+s)+s*2^s)*Int(t^(-s)*(s-2*s/t*sinh
(1/2*t)-1+cosh(1/2*t))/(exp(t)-1),t = 0 .. infinity);}{%
\[
\zeta (s)={\displaystyle \frac {\pi ^{s - 1}\,4^{s}\,
\mathrm{sin}({\displaystyle \pi\,s/\,2 } )}{2 - 2^{1 + s
} + s\,2^{s}}} \,{\displaystyle \int _{0}^{\infty }} 
{\displaystyle \frac {t^{ - s}\, \left (  \! s - {\displaystyle 
\frac {2\,s\,\mathrm{sinh}({\displaystyle t/2} )}{t}} 
 - 1 + \mathrm{cosh}({\displaystyle t/2} ) \;  \right ) 
}{e^{t} - 1}} \,dt \, . \hspace*{1cm} \sigma<2\,\,
\]
}
\label{IntRep1}
\end{equation}

Each of the terms in (\ref{IntRep1}) correspond to known integrals (\cite{G&R}, Eqs. (3.552.1) and (8.312.2)), and if the identification is taken to completion, (\ref{IntRep1}) reduces to the identity $\zeta(s)=\zeta(s)$ with the help of (\ref{functional}). Other useful representations are immediately available by evaluating the integrals in (\ref{IntRep1}) corresponding to pairs of terms in the large parenthesis, giving

\begin{equation}
\mapleinline{inert}{2d}{Zeta(s) =
1/4*Pi^s/(1-2^s)*2^s/cos(1/2*Pi*s)*(sin(Pi*s)*2^s/Pi*Int(t^(-s)*(-1+co
sh(1/2*t))/(exp(t)-1),t = 0 .. infinity)+1/GAMMA(s));}{%
\[
\zeta (s)={\displaystyle } \,{\displaystyle \frac {
\pi ^{s}\,2^{s-2}\,}{(1 - 2^{s})
\,\mathrm{cos}({\displaystyle \pi\,s/\,2} )}  \left(  \! {\displaystyle \frac {\mathrm{sin}(
\pi \,s)\,2^{s}}{\pi }} \,{\displaystyle \int _{0}^{\infty }} 
{\displaystyle \frac {t^{- s}\,( - 1 + \mathrm{cosh}(
{\displaystyle \frac {t}{2}} ))}{e^{t} - 1}} \,dt + 
{\displaystyle \frac {1}{\Gamma (s)}}  \!  \right) } \hspace*{1cm} \sigma<2
\]
}
\label{IntRep1a}
\end{equation}

and

\begin{equation}
\mapleinline{inert}{2d}{Zeta(s) =
-Pi^(-1+s)*sin(1/2*Pi*s)*(-Int(v^(-s-1)*(v-sinh(v))*exp(-v)/sinh(v),v
= 0 .. infinity)+Pi/sin(Pi*s)/GAMMA(1+s));}{%
\[
\zeta (s)= \pi ^{s - 1}\,\mathrm{sin}({\displaystyle 
\pi\,s/\,2} )\, \left(  \!  {\displaystyle \int _{0}^{
\infty }} {\displaystyle {v^{ - s - 1}\,(v/\mathrm{sinh}(v)  - 1 )} \,e^{- v}} \,dv - {\displaystyle 
\frac {\pi }{\mathrm{sin}(\pi \,s)\,\Gamma (1 + s)}}  \! 
 \right)\,. \hspace*{1cm} \sigma<2\,\,.
\]
}
\label{IntRep1b}
\end{equation}

%% file: s_eq_n.tex
Several of the results obtained in previous sections reduce to helpful representations when $s=n$. 

\subsection{Finite series representations}
From (\ref{SumInt2})  we find, for even values of the argument,

\begin{equation}
\mapleinline{inert}{2d}{Zeta(2*n) =
-1/2*1/(2^(-2*n)+n-1)*Sum((-1)^k*Pi^(2*k)*(2*n-2*k-1)*Zeta(2*n-2*k)/GA
MMA(2*k+2),k = 1 .. n);}{%
\[
\zeta (2\,n)= - {\displaystyle \frac {1}{2}} \,{\displaystyle 
\frac {{\displaystyle \sum _{k=1}^{n}} \,{\displaystyle \frac {(
-1)^{k}\,\pi ^{2\,k}\,(2\,n - 2\,k - 1)\,\zeta (2\,n - 2\,k)}{
\Gamma (2\,k + 2)}} }{2^{- 2\,n} + n - 1}} \,,
\]
}
\label{Zsum2n1}
\end{equation}

from (\ref{Zsum1a}) we find
\begin{equation}
\zeta(2n) = \frac{1}{2(2^{-2n}-1)}\sum _{k=1} ^{n} \frac {(-\pi^{2})^{k} (1-2^{1-2n+2k})\zeta(2n-2k)}{\Gamma(2k+1)}
\label{Zsum2n2}
\end{equation}

and, based on (\ref{TyHo}) we find (\cite{Milg1}) a similar, but independent result

\begin{equation}
\mapleinline{inert}{2d}{Z2nc :=
2*2^(-2*n)*Sum(-(1-2^(2*n-2*k-1))*(2*Pi)^(2*k)*zeta(2*n-2*k)/GAMMA(2*k
+2)*(-1)^k,k = 1 .. n)/(-1+2^(1-2*n));}{%
\[
\zeta(2\,n) = {\displaystyle 2 \, \frac {2^{- 2\,n}  }{ 1 - 2^{1 - 
2\,n}} \, 
 \! {\displaystyle \sum _{k=1}^{n}} \,( {\displaystyle \frac {(
1 - 2^{2\,n - 2\,k - 1})\,(2\,\pi )^{2\,k}\,\zeta (2\,n - 2\,
k)\,(-1)^{k}}{\Gamma (2\,k + 2)}} ) \! } 
\]
}
\label{Zsum2n3}
\end{equation}

three of an infinite number of such recursive relationships (\cite{Deeba}). The latter two may be identified with known results by reversing the sums, giving

\begin{equation}
\mapleinline{inert}{2d}{zeta(2*n) =
-1/2*(-Pi^2)^n/(2^(-2*n)+n-1)*Sum((-Pi^2)^(-k)*(2*k-1)*Zeta(2*k)/GAMMA
(2*n-2*k+2),k = 0 .. n-1);}{%
\[
\zeta (2\,n)= - {\displaystyle } \,{\displaystyle 
\frac {( - \pi ^{2})^{n}} {2
^{1- 2\,n} + 2\,n - 2}} \, \left( \, {\displaystyle \sum _{k=0}
^{n - 1}} \,{\displaystyle \frac { \,(2\,k
 - 1)\,\zeta (2\,k)}{( - \pi ^{2})^{\, k} \, \Gamma (2\,n - 2\,k + 2)}}  \!  \right) \,\, ,
\]
}
\end{equation}
thought to be new,
\begin{maplegroup}

\mapleresult
\begin{equation}
\mapleinline{inert}{2d}{zeta(2*n) =
1/2*(-Pi^2)^n/(2^(-2*n)-1)*Sum((-Pi^2)^(-k)*(1-2^(1-2*k))*Zeta(2*k)/GA
MMA(2*n-2*k+1),k = 0 .. n-1);}{%
\[
\zeta (2\,n)={\displaystyle } \,{\displaystyle 
\frac {( - \pi ^{2})^{n} }{2^{1- 2\,n} - 2}} \, \left( \,{\displaystyle \sum _{k=0}
^{n - 1}} \,{\displaystyle \frac { (1 - 2
^{1 - 2\,k})\,\zeta (2\,k)}{( - \pi ^{2})^{\,k} \Gamma (2\,n - 2\,k + 1)}}  \! 
 \right)
\]
}
\end{equation}

\end{maplegroup}
equivalent to \cite{SrivChoi2}, Eq.4.4(11), and
\begin{equation}
\mapleinline{inert}{2d}{zeta(2*n) =
2*(-Pi^2)^n/(1-2^(1-2*n))*Sum((1-2^(2*k-1))*(2*Pi)^(-2*k)*Zeta(2*k)*(-
1)^k/GAMMA(2*n-2*k+2),k = 0 .. n-1);}{%
\[
\zeta (2\,n)={\displaystyle \frac {2\,( - \pi ^{2})^{n}} {1 - 2^{1 - 2\,n}} \, \left( 
 \, {\displaystyle \sum _{k=0}^{n - 1}} \,{\displaystyle \frac {(
1 - 2^{2\,k - 1})\, \,\zeta (2\,k)\,(-1)^{k
}}{(2\,\pi )^{\,2\,k} \,\Gamma (2\,n - 2\,k + 2)}}  \!  \right) 
} 
\]
}
\end{equation}

equivalent to \cite{SrivChoi2}, Eq.4.4(9).

\subsection{Novel series representations}

From (\ref{ZetaSplit}), (\ref{zeta_{-}2}) and (\ref{EiSum_with_N}) with $N\geq 2, s=2$ find

\begin{equation}
\mapleinline{inert}{2d}{zeta(2) = Ly[2](exp(-2*omega))-4*Sum(arctan(omega/Pi/m)*Pi*m-omega,m
= 1 .. infinity)-2*omega*ln(-1+exp(2*omega))+3*omega^2+2*omega;}{%
\[
\zeta (2)=- 4\,\omega \left( 
 \! {\displaystyle \sum _{m=1}^{\infty }} \,({\displaystyle \frac{\pi \,m}{\omega}} \,\mathrm{arctan}(
{\displaystyle \frac {\omega }{\pi \,m}} )\, -1 )
 \!  \right) + \, {\mathrm{Li}_{2}}(e^{ - 2\,\omega })  - 2\,\omega \,\mathrm{log}( - 1 + e^{2\,\omega })
 + 3\,\omega ^{2} + 2\,\omega 
\]
}
\label{Zeta(2)}
\end{equation}

after identifying the hypergeometric function with $s=2$ in (\ref{zeta_{-}2}) from standard sources (e.g. \cite{Lebedev}, Chapter 9)\footnote{Compare the case N=0 with (\cite{Hansen}, Eq. (42.1.1))}. The series (\ref{Zeta(2)}) is absolutely convergent by Gauss' test for arbitrary values of $\omega$ and reduces to an identity in the limiting case $\omega \rightarrow 0$. Similarly for the case $s=3$ we find

\begin{align}
\mapleinline{inert}{2d}{zeta(3) =
2*Ly[2](exp(-2*omega))*omega+Ly[3](exp(-2*omega))+2*omega^2*Sum(-Pi^2*
ln(1+omega^2/Pi^2/m^2)*m^2/omega^2+1,m = 1 ..
infinity)-2*omega^2*ln(-1+exp(2*omega))-1/3*(-10*omega-3)*omega^2;}{%
\maplemultiline{
\zeta (3) 
= - 2\,\omega ^{2}\, 
 \! \left(\,{\displaystyle \sum _{m=1}^{\infty }} \, ( \,
{\displaystyle \frac {\pi ^{2}\,\, m^{2}}{\omega ^{2}}\, \mathrm{log}(1 + {\displaystyle 
\frac {\omega ^{2}}{\pi ^{2} \,m^{2} }} ) -\,1 }   
\,  )  \right ) + 2\, \omega  \,{\mathrm{Li}_{2}}(e^{ - 2\,\omega })+ {
\mathrm{Li}_{3}}(e^{ - 2\,\omega })  \\
\hspace*{1cm} - 2\,\omega ^{2}\,\mathrm{log}( - 1 + e^{2\,\omega }) + 
{\displaystyle \frac {(  10\,\omega  + 3)\,\omega ^{2}}{3}}  }
}
\label{Zeta(3)}
\end{align}

for $n=4$ the result is
\begin{equation}
\mapleinline{inert}{2d}{zeta(4) =
8/3*omega^3*Sum(Pi^3*arctan(omega/Pi/m)*m^3/omega^3-Pi^2*m^2/omega^2+1
/3,m = 1 ..
infinity)+2*Ly[2](exp(-2*omega))*omega^2+2*Ly[3](exp(-2*omega))*omega+
Ly[4](exp(-2*omega))-4/3*omega^3*log(-1+exp(2*omega))+7/3*omega^4+4/9*o
mega^3;}{%
\maplemultiline{
\zeta (4)={\displaystyle \frac {8}{3}} \,\omega ^{3}\, \left( 
 \, {\displaystyle \sum _{m=1}^{\infty }} \, (  \, 
{\displaystyle \frac {\pi ^{3}\,m^{3}}{\omega ^{3}}} \,\mathrm{arctan}({\displaystyle 
\frac {\omega }{\pi \,m}} )   - 
{\displaystyle \frac {\pi ^{2}\,m^{2}}{\omega ^{2}}}  + 
{\displaystyle \frac {1}{3}}\, )  \!  \right)  + 2\, \omega ^{2} \,{
\mathrm{Li}_{2}}(e^{ - 2\,\omega }) + 2\,\omega  \,{\mathrm{
Li}_{3}}(e^{ - 2\,\omega }) \\
\mbox{} + {\mathrm{Li}_{4}}(e^{ - 2\,\omega }) - 
{\displaystyle \frac {4}{3}} \,\omega ^{3}\,\mathrm{log}( - 1 + e
^{2\,\omega }) + {\displaystyle \frac {7\,\omega ^{4}}{3}}  + 
{\displaystyle \frac {4\,\omega ^{3}}{9}}  }
}
\label{Zeta(4)}
\end{equation}
and for $n=5$
 
\begin{equation}
\mapleinline{inert}{2d}{zeta(5) =
2/3*Sum(ln((Pi^2*m^2+omega^2)/Pi^2/m^2)*Pi^4*m^4/omega^4-Pi^2*m^2/omeg
a^2+1/2,m = 1 ..
infinity)*omega^4+4/3*Ly[2](exp(-2*omega))*omega^3+2*Ly[3](exp(-2*omeg
a))*omega^2+2*Ly[4](exp(-2*omega))*omega+Ly[5](exp(-2*omega))-2/3*omeg
a^4*ln(-1+exp(2*omega))+6/5*omega^5+1/6*omega^4;}{%
\maplemultiline{
\zeta (5)={\displaystyle \frac {2}{3}} \,\omega ^{4} \left(  \, 
{\displaystyle \sum _{m=1}^{\infty }} \; (  \, 
{\displaystyle {\frac{ \pi ^{4}\,m^{4}}{\omega^{4}}\, \mathrm{log}({1 +\, \displaystyle \frac { \omega ^{2}}{\pi ^{2}\,m^{2}}} )\, } }  - {\displaystyle \frac {\pi ^{2}\,m^{2}}{\omega ^{
2}}}  + {\displaystyle \frac {1}{2}}  \; )  \!  \right) \,
 + {\displaystyle \frac {4}{3}}\,\omega ^{3} \,{\mathrm{Li}_{2}}(e
^{ - 2\,\omega })  \\
\mbox{} + 2\,\omega ^{2}\,{\mathrm{Li}_{3}}(e^{ - 2\,\omega }) 
 + 2\,\omega  \,{\mathrm{Li}_{4}}(e^{ - 2\,\omega }) + {\mathrm{
Li}_{5}}(e^{ - 2\,\omega }) - {\displaystyle \frac {2}{3}} \,
\omega ^{4}\,\mathrm{log}( - 1 + e^{2\,\omega }) + 
{\displaystyle \frac {6\,\omega ^{5}}{5}}  
\mbox{} + {\displaystyle \frac {\omega ^{4}}{6}}\, . } 
}
\label{Zeta(5)}
\end{equation}

When applied to (\ref{EiSum_with_N}), the general case follows from the well-known identifications (\cite{Lebedev}, Eq.(9.8.5)) 
\begin{equation}
\mapleinline{inert}{2d}{hypergeom([1/2, 1],[3/2],-omega^2/Pi^2/m^2)/m^2 =
arctan(omega/Pi/m)*Pi/m/omega;}{%
\[
{\displaystyle  {\mathrm{_{2}F_{1}}({\displaystyle \frac {1
}{2}} , \,1 ; \,{\displaystyle \frac {3}{2}} ; \, - 
{\displaystyle \frac {\omega ^{2}}{\pi ^{2}\,m^{2}}} )}\, /\,{m^{2}}} =
{\displaystyle \frac {\pi \,\mathrm{arctan}({\displaystyle \frac {
\omega }{\pi \,m}} )\, }{m\,\omega }} 
\]
}
\label{Hn=0}
\end{equation}

and 

\begin{equation}
\mapleinline{inert}{2d}{hypergeom([1, 1],[2],-omega^2/Pi^2/m^2)/m^2 =
Pi^2*ln(1+omega^2/Pi^2/m^2)/omega^2;}{%
\[
{\displaystyle {\mathrm{_{2}F_{1}}(1, \,1 ; \,2; \, - 
{\displaystyle \frac {\omega ^{2}}{\pi ^{2}\,m^{2}}} )}\,/\,{m^{2}}} =
{\displaystyle \frac {\pi ^{2}}{\omega ^{2}}\,\mathrm{ln}(1 + {\displaystyle 
\frac {\omega ^{2}}{\pi ^{2}\,m^{2}}} )  }\; . 
\]
}
\label{Hn=1}
\end{equation}

By adding and subtracting the first $n-1$ terms in the series representation of the hypergeometric function in (\ref{zeta_{-}2}), we find, for the case $s=2n+2$

\begin{equation}
\mapleinline{inert}{2d}{hypergeom([1, n+1/2],[n+3/2],x)/m^2 =
(n+1/2)*(-omega^2/Pi^2/m^2)^(-n)*(2*arctan(omega/Pi/m)*Pi/m/omega-Sum(
(-omega^2/Pi^2/m^2)^k/(1/2+k),k = 0 .. n-1)/m^2);}{%
\maplemultiline{
{\displaystyle {\mathrm{_{2}F_{1}}(1, \,n + 
{\displaystyle \frac {1}{2}} ; \,n + {\displaystyle \frac {3}{2
}} , \,\, - 
{\displaystyle \frac {\omega ^{2}}{\pi ^{2}\,m^{2}}})}\,/\,{m^{2}}} = \\
(n + {\displaystyle \frac {1}{2}} )\,{( - {\displaystyle \frac {
\omega ^{2}}{\pi ^{2}\,m^{2}}} )}^{ - n}\, \left(  \! 
{\displaystyle \frac {2\,\pi }{m\,\omega } \,\mathrm{arctan}({\displaystyle \frac {
\omega }{\pi \,m}} ) }  - {\displaystyle \frac{1}{m^{2}} 
{{\displaystyle \sum _{k=0}^{n - 1}} \,{\displaystyle 
\frac {{{( - {\displaystyle \frac {\omega ^{2}}{\pi ^{2}\,m^{2}}} )}^{k} }
}{{ 1/2}  + k}} } }  \!  \right) 
 }
}
\label{Hs=2n}
\end{equation}

and for the case $s=2n+3$

\begin{equation}
\mapleinline{inert}{2d}{hypergeom([1, n+1],[n+2],-omega^2/Pi^2/m^2)/m^2 =
(n+1)*(-omega^2/Pi^2/m^2)^(-n)*(Pi^2*ln(1+omega^2/Pi^2/m^2)/omega^2-Su
m((-omega^2/Pi^2/m^2)^k/(k+1),k = 0 .. n-1)/m^2);}{%
\maplemultiline{
{\displaystyle {\mathrm{_{2}F_{1}}(1, \,n + 1 ; \,n + 2
; \, - {\displaystyle \frac {\omega ^{2}}{\pi ^{2}\,m^{2}}} )}\,/\,{m
^{2}}} = \\
(n + 1)\,{( - {\displaystyle \frac {\omega ^{2}}{\pi ^{2}\,m^{2}}
} )}^{ - n}\, \left(  \! {\displaystyle \frac {\pi ^{2}}{\omega ^{2}} \,
\mathrm{ln}(1 + {\displaystyle \frac {\omega ^{2}}{\pi ^{2}\,m^{2
}}} ) }  - {\displaystyle \frac{1}{m^2} {{\displaystyle 
\sum _{k=0}^{n - 1}} \,{\displaystyle \frac {{( - {\displaystyle 
\frac {\omega ^{2}}{\pi ^{2}\,m^{2}}} )}^{k}}{k + 1}} }} 
 \!  \right)  }
}
\end{equation}
for use in \eqref{ZetaSplit}, \eqref{zeta_{-}8} and \eqref{EiSum_with_N}.  I can find no direct record in the literature, with particular reference to works involving Barnes double Gamma function (e.g. \cite{ChoiSrivAdam} ), that such representations are known, although, as will be shown below, they are equivalent to known results.

\subsection{Connection with multiple Gamma function}\label{sec:MGamma}

From Section(\ref{sec:splitting}) with $\omega=1$ we have, after some simplification,

\begin{equation}
\mapleinline{inert}{2d}{zeta(2) = 5-4*Sum((-1)^m*Pi^(-2*m)*Zeta(2*m)/(2*m+1),m = 1 ..
infinity)+polylog(2,exp(-2))-2*ln(exp(2)-1);}{%
\[
\zeta (2)=5 - 4\, \left(  \! {\displaystyle \sum _{m=1}^{\infty }
} \,{\displaystyle \frac {(-1)^{m}\,\pi ^{( - 2\,m)}\,\zeta (2\,m
)}{2\,m + 1}}  \!  \right)  + \mathrm{Li_2}(e^{-2}) - 2
\,\mathrm{log}(e^{2} - 1).
\]
}
\end{equation}

Comparison of the above infinite series with \cite{Sriv}, Eq.3.4(464) using $z=i/\pi$ identifies
\begin{equation}
\mapleinline{inert}{2d}{ln(G(1+1/Pi*I)/G(1-I/Pi)) :=
1/4*I*(6+4*ln(2*Pi)+2*polylog(2,exp(-2))-4*ln(exp(2)-1)-1/3*Pi^2)/Pi =
.232662175652953*I;}{%
\maplemultiline{
\mathrm{log} \left(  \! {\displaystyle \frac {\mathrm{G}(1 + 
{\displaystyle i/\pi } )}{\mathrm{G}(1 - {\displaystyle 
i/\pi } )}}  \!  \right) = {\displaystyle {
{\displaystyle \frac {i}{4\,\pi}} \,(6 + 4\,\mathrm{log}(2\,\pi ) + 
2\,\mathrm{Li_2}(e^{-2}) - 4\,\mathrm{log}(e^{2} - 1) - 
{\displaystyle  {\pi ^{2}/3} } )}{}} = 
0.232662175652953\,i }
}
\end{equation}

where $\mathrm{G}$ is Barnes double Gamma function. From \cite{Sriv}, Section 1.3, unnumbered equation following 26, we also find
\begin{equation}
\mapleinline{inert}{2d}{ln(G(1+1/Pi*I)/G(1-I/Pi)) =
ln(2*Pi)/Pi*I-I/Pi+Sum(-2*I/Pi+k*ln((k+1/Pi*I)/(k-I/Pi)),k = 1 ..
infinity);}{%
\[
\mathrm{log} \left(  {\displaystyle \frac {\mathrm{G}(1 + 
{\displaystyle i/\pi} )}{\mathrm{G}(1 - {\displaystyle 
i/\pi } )}}  \!  \right) ={\displaystyle \frac {i\,\mathrm{
log}(2\,\pi )}{\pi }}  - {\displaystyle \frac {i}{\pi }}  + 
 \left(\,  {\displaystyle \sum _{k=1}^{\infty }} \, \left(  \! 
 - {\displaystyle \frac {2\,i}{\pi }}  + k\,\mathrm{ln} \left( 
 \! {\displaystyle \frac {k + {\displaystyle i/\pi } }{k
 - {\displaystyle i/\pi} }}  \!  \right)  \!  \right) 
 \!  \right) 
\]
}
\end{equation}
thereby recovering \eqref{Zeta(2)} after representing $arctan$ in \eqref{Zeta(2)} as a logarithm. In the case $s=3$, from \eqref{zeta_{-}7} and \eqref{EiSum_with_N} we have the series representation
\begin{equation}
\mapleinline{inert}{2d}{zeta(3) =
2*Li_2(exp-2)+Li_3(exp-2)-2*log(exp(2)-1)+13/3-2*Sum((-1)
^m*Zeta(2*m)*Pi^(-2*m)/(m+1),m = 1 .. infinity);}{%
\maplemultiline{
\zeta (3)=2\,\mathrm{Li_2}(e^{-2}) + \mathrm{Li_3}(e^{-2}) - 2\,\mathrm{log}(e^{2} - 1) + {\displaystyle 
\frac {13}{3}}  
\mbox{} - 2\, {\displaystyle \sum _{m=1}^{\infty }} \,
{\displaystyle \frac {(-1)^{m}\,\zeta (2\,m)\,\pi ^{ - 2\,m}}{m
 + 1}}   }
}
\end{equation}

from which similar representions for the logarithm of the triple Gamma function $\Gamma_3$ can be obtained by recourse to \cite{Sriv}, Eqs, 3.4(445) and (501).

%% file: critstrip.tex
In (\ref{Int_sinh^2a}) and (\ref{Int_sinh^2b}) with $0 \leq \sigma \leq 1$ let

\mapleresult
\begin{maplelatex}
\begin{equation}
\mapleinline{inert}{2d}{f = -Pi^{s-1}*sin(1/2*Pi*s)/(s-1);}{%
\[
f(s)= - {\displaystyle \frac {\pi ^{s - 1}\,\mathrm{sin}(
{\displaystyle \pi\,s/\,2} )}{s - 1}} 
\]
}
\end{equation}
\end{maplelatex}

\begin{maplelatex}
\begin{equation}
\mapleinline{inert}{2d}{g = 2^(s-1)/GAMMA(s+1);}{%
\[
g(s)={\displaystyle \frac {2^{s - 1}}{\Gamma (s + 1)}} 
\]
}
\end{equation}
\end{maplelatex}

\begin{maplelatex}
\begin{equation}
\mapleinline{inert}{2d}{h(t) = t^(1/2)*(1/(sinh(t)^2)-1/(t^2));}{%
\[
\mathrm{h}(t)=\sqrt{t}\,({\displaystyle \frac {1}{\mathrm{sinh} ^{2} t
}}  - {\displaystyle \frac {1}{t^{2}}} )
\]
}
\end{equation}
\end{maplelatex}

Then, after adding and subtracting, we find a composition of symmetric and anti-symmetric integrals about the critical line $\sigma=1/2$ :

\begin{maplegroup}
\mapleresult
\begin{maplelatex}
\begin{equation}
\mapleinline{inert}{2d}{Zeta(s) =
2*f*g/(g-f)*Int(h(t)*cos(rho*ln(t))*sinh((1/2-sigma)*ln(t)),t = 0 ..
infinity)+2*I*f*g/(f-g)*Int(h(t)*sin(rho*ln(t))*cosh((1/2-sigma)*ln(t)
),t = 0 .. infinity);}{%
\maplemultiline{
\zeta (s)={\displaystyle \frac {2\,f(s)\,g(s)}{g(s) - f(s)}} \,
{\displaystyle \int _{0}^{\infty }} \mathrm{h}(t)\,\mathrm{cos}(
\rho \,\mathrm{ln}(t))\,\mathrm{sinh}(({\displaystyle 1/2}  - \sigma )\,\mathrm{ln}(t))\,dt \\
\mbox{} + {\displaystyle \frac {2\,i\,f(s)\,g(s)}{f(s) - g(s)}} \,
{\displaystyle \int _{0}^{\infty }} \mathrm{h}(t)\,\mathrm{sin}(
\rho \,\mathrm{ln}(t))\,\mathrm{cosh}(({\displaystyle 1/2}  - \sigma )\,\mathrm{ln}(t))\,dt }
\label{subtracting}
}
\end{equation}
\end{maplelatex}
\end{maplegroup}

and simultaneously

\begin{maplegroup}
\mapleresult
\begin{maplelatex}
\begin{equation}
\mapleinline{inert}{2d}{Zeta(s) =
2*f*g/(f+g)*Int(h(t)*cos(rho*ln(t))*cosh((1/2-sigma)*ln(t)),t = 0 ..
infinity)+(-2*I*f*g/(f+g)*Int(h(t)*sin(rho*ln(t))*sinh((1/2-sigma)*ln(
t)),t = 0 .. infinity));}{%
\maplemultiline{
\zeta (s)={\displaystyle \frac {2\,f(s)\,g(s)}{f(s) + g(s)}} \,
{\displaystyle \int _{0}^{\infty }} \mathrm{h}(t)\,\mathrm{cos}(
\rho \,\mathrm{ln}(t))\,\mathrm{cosh}(({\displaystyle  1/2}  - \sigma )\,\mathrm{ln}(t))\,dt \\
\mbox{} - {\displaystyle \frac {2\,i\,f(s)\,g(s)}{f(s) + g(s)}} \,
{\displaystyle \int _{0}^{\infty }} \mathrm{h}(t)\,\mathrm{sin}(
\rho \,\mathrm{ln}(t))\,\mathrm{sinh}(({\displaystyle 1/2}  - \sigma )\,\mathrm{ln}(t))\,dt\;. }
\label{adding}
}
\end{equation}
\end{maplelatex}
\end{maplegroup}

Let 
\begin{equation}
\mathrm{C}(\rho )={\displaystyle \int _{0}^{\infty }} \mathrm{h}(t)\,\mathrm{cos}(
\rho \,\mathrm{ln}(t))\,dt
\end{equation}
and
\begin{equation}
\mathrm{S}(\rho )={\displaystyle \int _{0}^{\infty }} \mathrm{h}(t)\,\mathrm{sin}(
\rho \,\mathrm{ln}(t))\,dt \; .
\end{equation}

For $\zeta(s)=0$ on the critical line $\sigma=1/2$, both \eqref{subtracting} and \eqref{adding} can be simultaneously satisfied  when both $C(\rho)\neq 0$ and $S(\rho) \neq 0$ only if
\begin{equation}
{\displaystyle |f(1/2+i\,\rho) |^{2}\,=\,|g(1/2+i\,\rho)|^2}
\end{equation}
a true relationship that can be verified by direct computation. Set $\sigma=1/2$ in (\ref{subtracting}) and (\ref{adding}) and find

\begin{equation}
\zeta(1/2+i\,\rho) = \frac{2\,f(1/2+i\,\rho) \,g(1/2+i\,\rho)}{f(1/2+i\,\rho) \,+\,g(1/2+i\,\rho) }\,\mathrm{C}(\rho) = \frac{2\,i\,f(1/2+i\,\rho) \,g(1/2+i\,\rho) }{f(1/2+i\,\rho) -g(1/2+i\,\rho) }\,\mathrm{S}(\rho)
\label{Zeta_f=Zeta_g}
\end{equation}

Eq. (\ref{Zeta_f=Zeta_g}) can be used to obtain a relationship between the real and imaginary components of  $\zeta(1/2+i\rho)$. Let

\begin{equation}
\zeta(1/2+i\rho) = \zeta_{R}(1/2+i\,\rho) \,+\,i\,\zeta_{I}(1/2+i\,\rho).
\end{equation}

Then, by breaking $f(1/2+i\,\rho)$ and $g(1/2+i\,\rho)$ into their real and imaginary parts, (\ref{Zeta_f=Zeta_g}) reduces to
 
\begin{maplelatex}
\begin{equation}
\mapleinline{inert}{2d}{zeta[R](1/2+i*rho)/zeta[I](1/2+i*rho) =
(Cm*cos(rho*ln(2*Pi))-Cp*sin(rho*ln(2*Pi)))/(Cp*cos(rho*ln(2*Pi))+Cm*s
in(rho*ln(2*Pi))-Pi^(1/2));}{%
\[
{\displaystyle \frac {{\zeta _{R}}({\displaystyle 1/2} 
 + i\,\rho )}{{\zeta _{I}}({\displaystyle  1/2}  + i\,
\rho )}} ={\displaystyle \frac {\mathit{C_{m}(\rho)}\,\mathrm{cos}(\rho \,
\mathrm{ln}(2\,\pi )) - \mathit{C_{p} (\rho) }\,\mathrm{sin}(\rho \,\mathrm{
ln}(2\,\pi ))}{\mathit{C_{p} (\rho) }\,\mathrm{cos}(\rho \,\mathrm{ln}(2\,
\pi )) + \mathit{C_{m} (\rho) }\,\mathrm{sin}(\rho \,\mathrm{ln}(2\,\pi ))
 - \sqrt{\pi }}} 
\]
}
\label{crit_Zr/Zi}
\end{equation}
\end{maplelatex}

where

\begin{equation}
\mapleinline{inert}{2d}{C[p] =
Re(GAMMA(1/2+i*rho))*cosh(1/2*Pi*rho)+Im(GAMMA(1/2+i*rho))*sinh(1/2*Pi
*rho);}{%
\[
{C_{p} (\rho) }= \Gamma_{R} ({\displaystyle 1/2}  + i\,\rho ) \,
\mathrm{cosh}({\displaystyle {\pi \,\rho\,/2 } } ) + \Gamma_{I} ({\displaystyle 1/2}  + i\,\rho ) \,\mathrm{sinh}
({\displaystyle {\pi \,\rho /2 } } ) 
\]
}
\end{equation}
\begin{equation}
\mapleinline{inert}{2d}{C[p] =
Re(GAMMA(1/2+i*rho))*cosh(1/2*Pi*rho)+Im(GAMMA(1/2+i*rho))*sinh(1/2*Pi
*rho);}{%
\[
{C_{m} (\rho)}= \Gamma_{I} ({\displaystyle 1/2}  + i\,\rho ) \,
\mathrm{cosh}({\displaystyle {\pi \,\rho\,/2 } } ) - \Gamma_{R} ({\displaystyle 1/2}  + i\,\rho ) \,\mathrm{sinh}
({\displaystyle {\pi \,\rho /2 } } )
\]
}
\end{equation}
in terms of the real and imaginary parts of $\Gamma (1/2+i\,\rho)$. This is a known result (\cite{Milg3}) that depends only on the functional equation (\ref{functional})\footnote{See the Appendix for further discussion on this point.}. In \cite{Milg3} it was shown that the numerator and denominator on the right-hand side of (\ref{crit_Zr/Zi}) can never vanish simultaneously, so this result defines the half-zeros\footnote{ only one of $\zeta_{R}(1/2+i\,\rho)$ or $\zeta_{I}(1/2+i\,\rho)$ vanishes}  (and usually brackets the full-zeros) of $\zeta(1/2+I\,\rho)$ according to whether the numerator or denominator vanishes. This was studied in \cite{Milg3}. The full-zeros\footnote{$\zeta_{R}(1/2+i\,\rho_0) =\zeta_{I}(1/2+i\,\rho_0) =0$} of $\zeta(1/2+i\,\rho)$ occur when the following condition 
\begin{align}
& S(\rho_{0})=0,\,
\label{S=C=0}
\\ & C(\rho_{0}) = 0
\label{S=C=00}
\end{align}
is simultaneously satisfied for some value(s) of $\rho_{0}$. In that case, the defining equation (\ref{crit_Zr/Zi}) becomes

\begin{equation}
\mapleinline{inert}{2d}{zeta[R](1/2+i*rho)/zeta[I](1/2+i*rho) =
(Cm*cos(rho*ln(2*Pi))-Cp*sin(rho*ln(2*Pi)))/(Cp*cos(rho*ln(2*Pi))+Cm*s
in(rho*ln(2*Pi))-Pi^(1/2));}{%
\[
{\displaystyle \frac {  {\zeta _{R}}({\displaystyle 1/2} 
 + i\,\rho ) \, ^{\prime} }{ {\zeta _{I}}({\displaystyle  1/2}  + i\,
\rho ) \, ^{\prime} }} = - \, {\displaystyle \frac {  {\zeta \,^{\prime} }({\displaystyle 1/2} 
 + i\,\rho ) \, _{I}}{ {\zeta \,^{\prime}}({\displaystyle  1/2}  + i\,
\rho ) _{R}}} = {\displaystyle \frac {\mathit{C_{m}(\rho)}\,\mathrm{cos}(\rho \,
\mathrm{ln}(2\,\pi )) - \mathit{C_{p} (\rho) }\,\mathrm{sin}(\rho \,\mathrm{
ln}(2\,\pi ))}{\mathit{C_{p} (\rho) }\,\mathrm{cos}(\rho \,\mathrm{ln}(2\,
\pi )) + \mathit{C_{m} (\rho) }\,\mathrm{sin}(\rho \,\mathrm{ln}(2\,\pi ))
 - \sqrt{\pi }}} \;  
\]
}
\label{crit_Zr'/Zi'}
\end{equation}

where the derivatives are taken with respect to the variable $\rho$, and \eqref{crit_Zr'/Zi'} is evaluated at $\rho=\rho_{0}$. Many solutions to (\ref{S=C=0}) and \eqref{S=C=00} are known (\cite{Odlyzko}), so there is no doubt that these two functions can simultaneously vanish\footnote{A proof that this cannot occur in \eqref{subtracting} and \eqref{adding} when $\sigma\neq1/2$ is left as an exercise for the reader.} at certain values of $\rho$. \newline

A further interesting representation can be obtained inside the critical strip by setting $c=\sigma$ in \eqref{ContIntEta1} and separating $\eta(s)$ into its real and imaginary components

\begin{equation}
\mapleinline{inert}{2d}{Zeta(sigma+rho*I)[R] = 
1/2*sigma^(1-sigma)*Int((1+t^2)^(-1/2*sigma)*exp(rho*arctan(t))*(cos(O
mega(t))*sin(Pi*sigma)*cosh(Pi*sigma*t)-sin(Omega(t))*cos(Pi*sigma)*si
nh(Pi*sigma*t))/(cosh(Pi*sigma*t)^2-cos(Pi*sigma)^2),t = -infinity ..
infinity);}{%
\maplemultiline{
{\eta_{R} (\sigma  + \rho \,I) }= {\displaystyle \frac {1}{2}} 
\sigma ^{(1 - \sigma )}\\ \times {\displaystyle \int _{ - \infty }^{\infty 
}} { \displaystyle \frac {(1 + t^{2})^{ - \sigma/2}\,e
^{\rho \,\mathrm{arctan}(t)}\,(\mathrm{cos}(\Omega (t))\,
\mathrm{sin}(\pi \,\sigma )\,\mathrm{cosh}(\pi \,\sigma \,t) - 
\mathrm{sin}(\Omega (t))\,\mathrm{cos}(\pi \,\sigma )\,\mathrm{
sinh}(\pi \,\sigma \,t))}{\mathrm{cosh}^{2}(\pi \,\sigma \,t) - 
\mathrm{cos}^{2}(\pi \,\sigma )}}  
dt }
\label{eta_R}
}
\end{equation}
\begin{equation}
\mapleinline{inert}{2d}{Zeta(sigma+rho*I)[I] =
-1/2*sigma^(1-sigma)*Int((1+t^2)^(-1/2*sigma)*exp(rho*arctan(t))*(sin(
Omega(t))*sin(Pi*sigma)*cosh(Pi*sigma*t)+cos(Omega(t))*cos(Pi*sigma)*s
inh(Pi*sigma*t))/(cosh(Pi*sigma*t)^2-cos(Pi*sigma)^2),t = -infinity ..
infinity);}{%
\maplemultiline{
{\displaystyle \eta_{I}  (\sigma  + \rho \,I)}= - {\displaystyle \frac {1}{2}} 
\sigma ^{(1 - \sigma )}\\ \times {\displaystyle \int _{ - \infty }^{\infty 
}} {\displaystyle \frac {(1 + t^{2})^{- \sigma/2 }\,e
^{\rho \,\mathrm{arctan}(t)}\,(\mathrm{sin}(\Omega (t))\,
\mathrm{sin}(\pi \,\sigma )\,\mathrm{cosh}(\pi \,\sigma \,t) + 
\mathrm{cos}(\Omega (t))\,\mathrm{cos}(\pi \,\sigma )\,\mathrm{
sinh}(\pi \,\sigma \,t))}{\mathrm{cosh}^{2}(\pi \,\sigma \,t) - 
\mathrm{cos}^{2}(\pi \,\sigma )}}  
dt }
\label{eta_I}
}
\end{equation}
where
\begin{equation}
\mapleinline{inert}{2d}{Omega(t) = 1/2*rho*ln(sigma^2*(1+t^2))+sigma*arctan(t);}{%
\[
\Omega (t)={\displaystyle \frac {1}{2}} \,\rho \,\mathrm{log}(
\sigma ^{2}\,(1 + t^{2})) + \sigma \,\mathrm{arctan}(t)\,.
\]
}
\end{equation}

In analogy to \eqref{subtracting} and \eqref{adding}, the representations \eqref{eta_R} and \eqref{eta_I} include a term $\mathrm{sin}(\pi \,\sigma )$ that vanishes on the critical line $\sigma=1/2$, both reducing to the equivalent of \eqref{EtaC=half} in that case\footnote{It is suggested that such forms contain information about the distribution of zeros of $\eta(s)$ and hence $\zeta(s)$ inside the critical strip.}. The integrands in both representations vanish as 
$\exp(-\pi \,\sigma t)$ as $t\rightarrow\infty$ and are finite at $t=0$; most of the numerical contribution to the integrals comes from the region $0<t\cong\rho$. Numerical evaluation in the range $0 \leq t $ can be acclerated by writing 
\begin{equation}
\mathrm{e}^{\rho\,\mathrm{arctan}(t)}=\mathrm{e}^{\rho\,(\mathrm{arctan}(t)-\pi/2)},
\end{equation}
and extracting a factor $\mathrm{e}^{\rho \,\pi/2}$ outside the integral in that range. It remains to be shown how the remaining integral conspires to cancel this factor in the asymptotic range $\rho\rightarrow\infty$. Finally, these two representations are invariant under the transformations $\rho\rightarrow \,- \rho$ and $t \rightarrow -t$, demonstrating that $\eta(s)$ is self-conjugate, at least in the critical strip.

%% file: comments.tex
A number of related results that are not easily categorized, and which either do not appear in the reference literature or which appear incorrectly should be noted.\newline 

Comparison of (\ref{Sum1_s}) and (\ref{PsiDiff0}) suggests the identification $v=i$. Such sums and integrals have been studied (\cite{Adam1}, \cite{ChoiSrivAdam}) for the case $s=-n$ in terms of Barnes double Gamma function $G(z)$ and Clausen's function defined by
\begin{equation}
{\displaystyle \mathrm{Cl_{2}}(z) =\sum_{k=1}^{\infty} \frac{ \mathrm{sin}(z\,k)}{k^2}},
\end{equation}

with the (tacit) restriction that $z \in \mathfrak{R}$. That is not the case here. Thus straightforward application of that analysis will sometimes lead to incorrect results. For example \cite{ChoiSrivAdam}, Eq.(4.3) gives (with $n=1$)

\begin{equation}
\pi\,\int_{0}^{t} v^{n} \, {\displaystyle \mathrm{cot}}(\pi v) dv = t\,\, {\displaystyle \mathrm{log} }(2\, \pi)\,+\,\mathrm{log} \left( \frac{G(1-t)}{G(1+t)} \right)
\label{v^n cot}
\end{equation}

further (\cite{Adam1}, Eq.(13)) identifying
\begin{equation}
\mathrm{log} \left( \frac{G(1-t)}{G(1+t)} \right)= t\,\mathrm{log} \left( \frac { \mathrm{sin}(\pi \,t)}{\pi} \right )  +\frac{1}{2\,\pi} \,\mathrm{Cl_{2}}(2 \pi t) \hspace*{1cm} 0<t<1 \:.
\end{equation}

To allow for the possibility that $t\in \mathfrak{C}$, this should instead read
\begin{equation}
\mathrm{log} \left( \frac{G(1-t)}{G(1+t)} \right)= t\,\mathrm{log} \left( \frac { \mathrm{sin}(\pi \,t)}{\pi} \right ) + \frac{i\pi}{2}\mathit{B_{2}}(t)+ \frac{\mathrm{e}^{-i\pi/2}}{2\pi}\mathrm{Li_{2}}(\mathrm{e}^{2\pi i t})
\label{G_with_Li} 
\end{equation}
valid for all $t$ in terms of Bernoulli polynomials $B_{n}(t)$ and polylogarithm functions. However, similar results $\mathit{are}$ given in the equivalent, general form in \cite{SG&A}, Eq. (5.12). Interestingly, both the Maple and Mathematica computer codes manage to evaluate\footnote{Neither \cite{Cv&Sriv}, nor references contained therein, both of which study such integrals, list the computer results, which are likely based on \cite{SG&A}, Eq. (5.12).} the integral in (\ref{v^n cot}) in terms similar to (\ref{G_with_Li}) for a variety of values of $n$ provided that the user restricts $|t|<1$, although the right-hand side provides the analytic continuation of the left when that restriction is lifted. This has application in the evaluation of $\zeta(3)$ according to the method of Section \eqref{sec:MGamma}, where, using the methods of \cite{Sriv}, Chapter 3, we identify 

\begin{equation}
{\displaystyle \sum _{k=1}^{\infty }} \,{\displaystyle \frac {
\zeta (2\,k)\,t^{(2\,k)}}{k + 1}} ={\displaystyle \frac {1}{2}} 
 - {\displaystyle \frac {\pi }{t^{2}}} \,{\displaystyle \int _{0
}^{t}} v^{2}\,\mathrm{cot}(\pi \,v)\,dv,
\label{Intv^2}
\end{equation}

a generalization of integrals studied in \cite{Cv&Sriv}. From the computer evaluation of the integral, \eqref{Intv^2} eventually yields

\begin{equation}
\mapleinline{inert}{2d}{Sum(Zeta(2*k)*t^(2*k)/(k+1),k = 1 .. infinity) =
1/3*I*Pi*t+1/2-ln(1-exp(2*I*t*Pi))+polylog(2,exp(2*I*t*Pi))/Pi/t*I+(-1
/2*polylog(3,exp(2*I*t*Pi))+1/2*Zeta(3))/Pi^2/t^2;}{%
\maplemultiline{
{\displaystyle \sum _{k=1}^{\infty }} \,{\displaystyle \frac {
\zeta (2\,k)\,t^{2\,k}}{k + 1}} ={\displaystyle \frac {1}{3}} 
\,i\,\pi \,t + {\displaystyle \frac {1}{2}}  - \mathrm{log}(1 - e
^{2\,i\,\pi\, t }) + i\,{\displaystyle \frac {\mathrm{Li_2}( 
\,e^{2\,i\,\pi\,t})}{\pi \,t}} 
\mbox{} - {\displaystyle \frac { {\displaystyle } 
\,\mathrm{Li_3}(e^{2\,i\,\pi\,t }) - {\displaystyle 
} \,\zeta (3)}{2\,\pi ^{2}\,t^{2}}}  }
}
\end{equation}
as a new addition to the lists found in \cite{Sriv}, Section (3.4) and \cite{SrivChoi2}, Sections (3.3) and (3.4), and a generalization of \eqref{zeta_{-}7} when $s=3$.

%% file: appendix.tex
In \cite{Milg3}, an attempt was made to show that the functional equation \eqref{functional} and a minimal number of analytic properties of $\zeta(s)$ might suffice to show that $\zeta(s) \neq 0$ when $\sigma \neq 1/2$. This was motivated by a comment attributed to H. Hamburger (\cite{Ivic}, page 50) to the effect that the functional equation for $\zeta(s)$ perhaps characterizes it completely\footnote{An unequivocal equivalent claim can also be found in \cite{Encyclop}. As demonstrated by the counterexample given here, and in \cite{Gauthier}, this claim is false.}. The purported demonstration contradicts a result of Gauthier and Zeron (\cite{Gauthier}), who show that infinitely many functions exist, arbitrarily close to $\zeta(s)$, that satisfy \eqref{functional} and the same analytic properties that are invoked in \cite{Milg3}, but which vanish when $\sigma \neq 1/2$.
\footnote{To quote Prof. Gauthier: "We prove the existence of functions approximating the Riemann zeta-function, satisfying the Riemann functional equation and violating the analogue of the Riemann hypothesis."(private communication)}\newline 

In fact, an example of such a function is easily found. Let 
\begin{equation}
\nu(s) = C\,\mathrm{sin}(\pi\,(s-s_0))\,\mathrm{sin}(\pi\,(s+s_0))\,\mathrm{sin}(\pi\,(s-s_0^*)\,\mathrm{sin}(\pi\,(s+s_0^*)) 
\end{equation}
where $s_0^*$ is the complex conjugate of $s_0$, $s_0=\sigma_0+i\rho_0$ is arbitrary complex, and $C$, given by 
\begin{equation}
C = {\displaystyle \frac {1}{(\mathrm{cosh}^{2}(\pi \,\rho_0) - 
\mathrm{cos}^{2}(\pi \,\sigma_0))^{2}}} 
\end{equation}

is a constant chosen such that $\nu(1)=1$. 
Then 
\begin{equation}
\zeta_{\nu}(s) \equiv \nu(s)\,\zeta(s) 
\end{equation}
satisfies analyticity conditions (1-3) of \cite{Gauthier} (which include the functional equation \eqref{functional}). Additionally, $\zeta_{\nu}(s)$ has zeros at the zeros of $\zeta(s)$ as well as complex zeros at $s=\pm\,(n+s_0)$ and $s=\pm\,(n+s_0^*)$ and a pole at $s=1$ with residue equal to unity, where $n=0,\pm 1,\pm 2,\dots $. \newline

Thus the last few paragraphs of Section 2 of \cite{Milg3} beginning with the words "There are two cases where (2.9) fails... " cannot be  valid. However, the remainder of  that paper, particularly those sections dealing with the properties of $\zeta(1/2+i\rho)$ as referenced here in Section \ref{sec:crits} are valid and correct. I am grateful to Prof. Gauthier for discussions on this topic.